\documentclass[poms,nonblindrev]{poms1} % Use this for blind submissions

\OneAndAHalfSpacedXI % % default line spacing !!!!!!!!!!!!!!!!!!!!!!!
\usepackage{graphicx}
\usepackage{natbib}
\usepackage{multirow}
\usepackage{endnotes}
\usepackage{bm}
\usepackage{mathrsfs}
\usepackage[boxed]{algorithm2e}
\usepackage{subfig}
\usepackage[colorlinks,citecolor=green]{hyperref}
\usepackage{easyReview}

\bibpunct[, ]{(}{)}{,}{a}{}{,}%
\def\bibfont{\small}%
\TheoremsNumberedThrough     % Preferred (Theorem 1, Lemma 1, Theorem 2)
\ECRepeatTheorems

\EquationsNumberedThrough    % Default: (1), (2), ...
\begin{document}
	%%%%%%%%%%%%%%%%
	
	%default values will be used.
	%\setcounter{page}{1}
	%\DOI{10.1287/xxxx.0000.0000}%
	
	%\RUNAUTHOR{Slippery and Arinella}
	
	
	%\RUNTITLE{Tis a Butter Place}
	
	% \TITLE{Bundling Information Goods of Decreasing Value}
	\TITLE{Explainable Data-driven Share-of-choice Product Line Design Optimization}
	\RUNTITLE{An Explainable PLD Model}
	%should be entered in ONE field, separated by a comma.
	%\EMAIL field can be repeated if more than one author
\ARTICLEAUTHORS{%
	\AUTHOR{Maoqi Liu\textsuperscript{1}, Xun Zhang\textsuperscript{2}, Hailei Gong\textsuperscript{3}, Changchun Liu\textsuperscript{4}}
	\AFF{\textsuperscript{1}School of Management, Shandong University, Jinan, China, 250061} %, \URL{}}
\AFF{\textsuperscript{2}College of Business,  Southern University of Science and Technology, Shenzhen, China, 518055}
\AFF{\textsuperscript{3} Department of Industrial Engineering, Tsinghua University, Beijing, China, 100084}
\AFF{\textsuperscript{4} School of Management, Xi'an Jiaotong University, Xi'an, China, 710049}
\AFF{\EMAIL{liumq@sdu.edu.cn, xunzhang@u.nus.edu, gonghl18@mails.tsinghua.edu.cn, liuchangchun@xitu.edu.cn}}
} % end of the block

\ABSTRACT{%
	
% PLD
%Product line design (PLD) focuses on determining the specifications of a set of products to meet diverse customer preferences.
% SOC Problem and Utility
The share-of-choice (SOC) problem is a widely studied problem for product line design (PLD) where representative customers are sampled from a target population, and the goal is to maximize the percentage of customers who choose the offered products over outside options over the sample. The individual choices are captured by the utility maximization framework. 
% Common Approach
A significant challenge is that these utilities are not directly observable and must be estimated from other primitive data. Conjoint analysis is a commonly applied technique for generating such data, where sampled customers rate, rank, or choose between different product alternatives. With the responses, various methods, such as the hierarchical Bayesian method and polyhedral estimation, are employed to estimate the utilities.
% Challenge
However, this ``estimate-then-optimize'' procedure disconnects the decision-making process from the primitive data and thus is unable to attribute the PLD decision to the survey data unresolved.
% Our work
%% Main Frame
To fill the gap, we integrates the polyhedral estimation method proposed by \cite{Toubia2003Fast,toubia2004polyhedral}, which constructs a polyhedron set containing all utility consistent with the responses, into the PLD optimization. Specifically, we proposed a robust model that maximizes the share-of-choice calculated by the worst-case utility over the polyhedron. The model brings threefold explainability to the PLD problem.
%%%1.
First, it enables us to explore the influence of both estimation error and the number of sampled customers through an out-of-sample performance guarantee of the proposed model. 
%%% 2. 
Secondly, the linearly tractable equivalent reformulation constructed by the duality allows us to attribute the optimal product line to the survey data through the dual variables. 
%%% 3. 
Third, inspired by the fact that the survey questions are associated with columns in the dual reformulation, we propose a column-generation approach to identify new questions that most effectively improve the model's objective. 
%%% A common summary.
Through extensive numerical results, we show the effectiveness of our proposed model and methods.

}%

%\KEYWORDS{deterministic inventory theory; infinite linear programming duality;
%existence of optimal policies; semi-Markov decision process; cyclic schedule}

%\HISTORY{This paper was first submitted on April 12, 1922 and has been with the authors for 83 years for 65 revisions.}

\maketitle
%%%%%%%%%%%%%%%%%%%%%%%%%%%%%%%%%%%%%%%%%%%%%%%%%%%%%%%%%%%%%%%%%%%%%%

% Samples of sectioning (and labeling) in OPRE
% NOTE: (1) \section and \subsection do NOT end with a period
%(2) \subsubsection and lower need end punctuation
%(3) capitalization is as shown (title style).
%
%\section{Introduction.}\label{intro} %%1.
%\subsection{Duality and the Classical EOQ Problem.}\label{class-EOQ} %% 1.1.
%\subsection{Outline.}\label{outline1} %% 1.2.
%\subsubsection{Cyclic Schedules for the General Deterministic SMDP.}
%\label{cyclic-schedules} %% 1.2.1
%\section{Problem Description.}\label{problemdescription} %% 2.

% Text of your paper here

\section{Introduction}\label{sec_Intro}

Product line design (PLD) is a widely adopted marketing practice aimed at accommodating diverse customer preferences. It has received substantial attention from the fields of operations research and marketing.The share-of-choice (SOC) problem is a focal point of particular prominence in both early \citep{10.2307/183706, Kohli&Sukumar1990, Balakrishnan1996Genetic} and contemporary literature \citep{Wang2012A, wang2022robust, liu2023share}.

The problem is established based on a sample of representative customers. The target of the problem is to find a product line maximizing the percentage of customers, who choose products from the offered line over the outside option, in this sample. In the problem, the customer choice is modeled by first-choice model. Specifically, for each individual customer, a utility is associated to each product and the outside option and  the customer the alternative with highest utility. It is challenging to get the utilities for two reasons.
\begin{itemize}
	\item It is almost impossible to get the utilities for all candidate product alternative. There are always many product alternatives to offer and therefore, it requires extensive efforts to assign the utility value for all of them. Previous literature tackled this issue by decomposing the products into attributes and assuming the utility of a product is the sum of the partworth utility associated to the performance levels selected for its attributes. Accordingly, it only need the partworths for the levels, whose number is much less than that of products, to calculate the utility of the products.
	\item The utilities are not directly measurable. That is, the utilities need to be calibrated from other primitive data.  \textit{Conjoint analysis (CA)} is a commonly applied set of techniques used to collect customer survey data and estimate the partworth. In CA, sampled customers are presented with a variety of products characterized by different attributes and are asked to rate, rank, or choose their preferred options from the presented alternatives \citep{Rao2014Applied}. Several methods have been proposed to infer partworth utilities, including the polyhedral method \citep{Toubia2003Fast, toubia2004polyhedral} and the hierarchical Bayesian method \citep{allenby1998marketing, Gilbride2004}. These methods can provide point estimations of utilities, which are subsequently input into the SOC model.
\end{itemize}

While the ``estimate-then-optimize" scheme offers a feasible approach to implement the SOC model, it manually bring up a barrier between the primitive data to the optimization procedure. Some information in these primitive data may lose during the estimation. As a result, the recommended product line derived from the SOC model lacks clear attribution to the primitive data. However, transparent reasoning connecting model decisions to the underlying data is important in the analytic, which enhances credibility and facilitates informed decision-making by stakeholders \citep{chen2022robust,de2023explainable}.

In this paper, we propose a end-to-end model, which directly input the CA survey data and out the recommended product line.  To achieve this, we use the the polyhedral method proposed by \cite{Toubia2003Fast,toubia2004polyhedral} to construct a polyhedral feasible set containing utilities consistent with each customer's responses and then propose a robust optimization approach to maximize the SOC with the worst-case utility for each customer over these polyhedral sets. Through duality, we present an equivalent tractable reformulation of the proposed model as a mixed integer linear programming, in which each CA survey question and response correspond to the coefficients of a dual variable. Expanding on this structure, this variables explains to what extent the optimal product line relies on the CA survey questions and responses. Moreover, inspired by the structure, we introduce an adaptive question design approach based on column generation, which aims to identify new survey questions that can most effectively improve the objective function. Extensive numerical experiments are conducted to validate the effectiveness of our proposed model and the adaptive design approaches.

To the best of our knowledge, our research is the first model bringing the primitive data into the product line optimization, which allows us better understand how the primitive data influence the subsequent optimization. There are threefold contributions of our study.
\begin{itemize}
	\item First, we provide a performance guarantee of the proposed model as a function of both the size of customer size and the estimation accuracy of the partworth. The guarantee uncovers the influence of the estimation accuracy on the PLD optimization, which is unable to be captured by the previous models.
	\item Second, the dual variables provide a way to attribute the optimal product line to primitive data. This enables a better understanding of how the decision is made according to the data collected and consequently, improve the credibility of the model to the decision maker.
	
	\item Third, we propose an optimization-directed adaptive design framework for CA survey. Adaptive design of CA survey is always based on the idea of improving the estimation accuracy as soon as possible and neglects its effect on the subsequent optimization.  Through extensive numerical experiments using synthetic data, we show the importance of aligning the CA design to the optimization. We show the superior performance of the proposed model and adaptive design by comparing it to the classic orthogonal designs and the Fast Polyhedral Adaptive Conjoint Analysis (FPACA) approach \citep{Toubia2003Fast,toubia2004polyhedral}.
	%\item Our model offers a robust solution for scenarios where the number of CA survey questions is limited relative to the candidate attribute levels and the sample size of customers is constrained. By leveraging worst-case utility considerations over polyhedral sets, our approach ensures robust decision-making with limited CA survey questions. Furthermore, we establish connections between our proposed model and recent advances in robust SOC maximization \citep{liu2023share} to illustrate its ability to offer robustness against uncertainty caused by limited number of sampled customers.
\end{itemize}

The remainder of this paper is organized as follows. Section \ref{sec_LR} reviews previous literature. Section \ref{sec_MF} provide preliminaries on the SOC problems, conjoint analysis and polyhedral estimation. Sections \ref{sec:end-to-end} and \ref{sec:ACA} present and analyze the proposed model and ACA framework, respectively. Section \ref{sec_NE} presents numerical experiments to validate the effectiveness of the proposed model and ACA framework.  Finally, conclusions and future research directions are outlined in Section \ref{sec_Con}. We furnish all the proofs and other technical supplements in the appendices for the briefness of reading.

\section{Literature Review}\label{sec_LR}

Our research intersects with two prominent streams of literature: Product Line Design (PLD) and the Adaptive Design of Conjoint Analysis (ACA).  In this section, we review the most relevant studies to showcase how our research addresses gaps in existing literature.

\subsection{Product line design based on First-choice Model}

The PLD problem is one of the central topic at the intersection of marketing and operations research. The goal of PLD is to find a product line maximizing specific market goals, such as profit, revenue, and market share. A basis to calculate all the objectives is the choice model, which describes customer choices when confronted with the offered product line. In this section, we focus on the first-choice model, which is one of the most widely applied one in the PLD literature, from the seminal works \citep{10.2307/183706,Kohli&Sukumar1990} to the extensive follow-up researches \citep{Shi2001An,camm2006conjoint,wang2009branch,Wang2012A,wang2022robust,liu2023share}. 

The first choice model characterize the choice of a sample of customers following the utility maximization manner. For each individual in the sample, each candidate product is associated with a deterministic utility and assumes that customers will select the option with the highest utility including both offered products and outside option \citep{bertsimas2016robust}. Note that although the utility is deterministic for each individual, the model is a sample average approximation (SAA) approach at the population level and provide a very flexible way to capture customer preference heterogeneity \citep{liu2023share}.  In early works \citep{10.2307/183706,Kohli&Sukumar1990,Shi2001An,Belloni2008Optimizing}, these utilities are assumed to be accurate. However, it is challenging to get the exact values of utilities is because they are not directly observable. In practice, they are estimated from some primitive data, usually the CA survey data, and therefore  are subjected to the estimation error. Recent studies have increasingly focused on developing robust models against this uncertainty. For example, \cite{Wang2012A} enhanced the SOC by considering the sample utilities within a specified budget uncertainty set. Furthermore, \cite{wang2022robust} introduced robust scenario-based robust model that characterize the uncertainty by multiple utility scenarios and optimizes the worst-case objectives  over scenarios. \cite{liu2023share} utilized Wasserstein distributionally robust optimization to maximize the worst-case objective across an uncertainty set comprising all utility distribution closed enough to the sample one with respect to the Wasserstein distance.

However, the uncertainty considered in these models are technically added and are not directly linked to its origin, i.e., the primitive data and the estimation procedure. While \cite{wang2022robust} suggested constructing scenarios through uniform bootstrapping from the polyhedrons constructed by the CA survey, the numerical nature of the approach limited the exploration of how the  data and estimation procedure influence the model’s solutions. In this work, we fill the gap by integrating the polyhedral estimation method proposed by \cite{Toubia2003Fast} into the optimization, which enable our model to nest the primitive CA survey data into optimization. Leveraging on the desirable geometric explainability of the polyhedral method and dual analysis, we obtain a tractable mixed integer linear formulation that can attribute of the solutions to the primitive CA data. Furthermore, our analysis on the performance guarantee explicitly shows how the estimation error influence the out-of-sample performance of the proposed model. To the best of our knowledge, our guarantee is the first one to capture this influence.

Other than the PLD, \cite{chen2023model} applied a similar idea to a data-driven pricing problem where a firm offers multiple pre-specified products to customers with varied prices. Transaction data recording the prices customer observed and the choice made, i.e., a similar format as choice-based CA survey, are applied to construct the polyhedral feasible set of utilities and worst-case revenue over the polyhedrons is maximized. The model is different to ours in three aspects. First, their pricing problem is with fixed product line and do not consider the product features. Accordingly, the model is not applicable to the case with varied product line. Second, their reformulation is based on the disjunctive programs and therefore, can not show the attribution of the product line to the primitive data through the dual analysis done in this paper. Third, their model only use the historical data statically while we leverage the structure of our model to design an adaptive approach that actively pursue new data.

\subsection{Adaptive Conjoint Analysis}

The PLD problems based on the first-choice model requires point estimations of partworth utilities at individual levels as input, which is usually an output of CA.  Careful design of the CA survey is a prerequisite of the success of the whole CA process \citep{netzer2008beyond}.  The CA design can be divided into two schemes, i.e., the static and adaptive ones. In this paper, we focus on the adaptive CA. As such, we neglect the static design approach and refer the readers to \cite{liu2015construction} for a thorough review. 

In ACA, the questions are designed sequentially during the survey and a feedback loop is built to use the responses to design next questions. \cite{Toubia2003Fast} is one of the pioneer researches on designing ACA approach. In the paper, they focused on the metric comparison CA, where the response is the exact utility gap between the two compared product profiles. They applied the CA survey data to construct polyhedral feasible sets of all partworth utilities consistent with the responses. Their ACA method designs the next-round question most parallel to the longest vectors in the polyhedrons so that the size of the feasible is reduced as fast as possible.   \cite{toubia2004polyhedral} extended the same idea to the choice-based CA. To involve the possible choice error during the survey, \cite{Toubia2007Probabilistic} proposed a probabilistic extension of \cite{toubia2004polyhedral}, which allows the choice error happening with certain probability. \cite{abernethy2007eliciting} combined the utility balance criterion, namely, the compared products are likely to be equal like by the respondent from the prior perspective, to the original polyhedral-based ACA approach. \cite{Bertsimas2013Learning} dealt with the respondent error to allow a certain number of the constraints to flip. \cite{saure2019ellipsoidal} applied an ellipsoidal credibility region to replace the polyhedron feasible set. They designed the new question to the expected the root of determinant of the covariance matrix, which will reduce the credibility region as fast as possible. These work differs from our ACA framework from two aspects. In contrast to our ACA framework, the main criterion of the approaches is to reduce the uncertainty in the partworth estimation while do not consider the objective of the subsequent optimization. Instead, we propose a goal-directed approach aiming at improving the objective of the subsequent optimization as soon as possible. 

A recent paper, \cite{JooOptimal}, also proposed a goal-directed ACA aiming at selecting a new question that compares the product profiles with high potentiality to improve the market share. The difference between this framework and ours is twofold. From estimation aspect, their approach only use the data on an aggregated level. Specifically, they use the share of an alternative to calibrate the parameters of a multi-nomial logit (MNL) choice model, which is the mean utility of the alternative over the population.  In contrast, our model take advantages of the data on both levels through a conditional probability on the individual choices. From the estimation perspective, their framework is an estimate-then-optimize approach, i.e., the utility is first estimated and then the optimal product is selected according to the estimated utility.  In contrast, ours utilize polyhedral method to direct take the primitive survey questions and response into the optimization, which leverage explainability to our model. 

The ACA framework driven by the optimization objective shares similar idea to the dynamic pricing \citep{besbes2009dynamic,farias2010dynamic,handel2015robust} or assortment planning \citep{caro2007dynamic,rusmevichientong2010dynamic,saure2013optimal,agrawal2019mnl}, where the pricing and assortment are dynamically adjusted to learn the customer preference and obtain good profit in the same time over the selling season. The stream of literature differs from ours from two aspects. From setting perspective, the offered products are considered as fixed in pricing and assortment planning and most of them do not consider the detail specifications of the products.  Although some literature, such as \cite{amin2014repeated,cohen2020feature}, considered product feature in the pricing problem, the features are exogenous covariates instead of decision in the optimization. From data perspective, the models only use the data to estimate the parameters of choice model at a population level instead of characterizing the individual partworth. Our framework utilizes the data at both level to better take advantages of the hierarchical structure of the CA data.

\section{Preliminary} \label{sec_MF}
\subsection{The Classic Model}
 In this section, we present the conventional setting and notations of the SOC problem in the previous literature \citep{Kohli&Sukumar1990,Shi2001An,camm2006conjoint,wang2009branch,Wang2012A,wang2022robust,liu2023share}. We consider a product line design problem where the decision-maker aims at launching $M$ products to the market. Each product is decomposed into $A$ attributes. There are several candidate levels for each attribute. Let $\mathcal{L}_a$ be the set of candidate levels for attribute $a\in [A]$ and $L=|\cup_{a\in{A}}\mathcal{L}_a|$ denote the number of all different levels. The design of the product line is then simplified as selecting one level for each attribute of each product in the line. Binary variables $x_{ml}$ represent the level selections where $x_{ml}=1$ if level $l$ is chosen for product $m$, and 0 otherwise. Accordingly, product $m\in [M]$ can be represented by a $L$-dimensional vector $\boldsymbol{x}_m=\{x_{1m},\cdots,x_{ml}\}$ and the designed product line design can be denoted by the matrix $\boldsymbol{X}= \{\boldsymbol{x}_1,\cdots,\boldsymbol{x}_{M}\}^T$. The SOC problem aims to find the optimal product line from a feasible set $\mathcal{X}$ so that the SOC is maximized. Here, we focus on a basic format where $\mathcal{X}=\{\boldsymbol{X}|\sum_{l\in\mathcal{L}_a}x_{ml}=1,\forall m\in[M],a\in[A]\}$, which requires that one and only one level is selected for each attribute of each product.
 
To get a proxy of the SOC over customer population, sample average approximation (SAA) approach and a utility maximization framework are applied. Suppose that we have $N$ sampled customers and for customer $n$, a partworth utility $\hat{u}_{nl}$, $n\in[N]$, $l\in[L]$ is associated to level $l$ to represent the customer's preference to the level. We assume the partworth utilities are additive, i.e., the utility of a product $\boldsymbol{x}$ for customer $n$ is $\boldsymbol{u}_n^T\boldsymbol{x}$ and the customer will choose the alternative  with highest utility over the union of offered products and outside option.  Without losing generality, we set the utility of the outside option by 0. Then, the SOC is defined as the percentage of customers, for whom at least one product in the line delivers a non-negative utility.  With above setting and notations, the classic SOC model \citep{Kohli&Sukumar1990,Shi2001An,camm2006conjoint,wang2009branch} is formally formulated as follows:
\begin{align}
	\mbox{[SOC-D]} \quad\max \limits_{\boldsymbol{X}\in \mathcal{X}}\frac{1}{N} \sum_{n \in [N]} \boldsymbol{1}\{\max_{m\in[M]} \boldsymbol{\hat{u}}^T_n\boldsymbol{x}_m\geq 0\}, \label{ClassicObj}
\end{align}

Equivalently, model SOC-D can be transformed into the following mixed integer linear programming.
\begin{eqnarray}
	\mbox{(SOC-D)}& \max\limits_{\boldsymbol{X},\boldsymbol{y}}&\frac{1}{N} \sum_{n \in [N]} y_n\\ 
	&s.t.& \boldsymbol{\hat{u}}^T_n\boldsymbol{x}_m+C(1-z_{nm})\geq 0,\forall n\in[N],\\ 
	&& \sum_{m\in[M]} z_{nm}\geq y_n,\forall n\in[N]\\
	&&  \sum_{l\in\mathcal{L}_a}x_{ml}=1,\forall a\in[A],m\in[M],\\
	&& x_{ml},y_n\in\{0,1\},\forall l\in[L],m\in[M],n\in[N].
\end{eqnarray}
where $C$ is a large enough constant number; $z_{nm}$ is a binary variable indicating whether product $m$ deliver a non-negative higher utility than the hurdle one; while $y_n$ implies whether there are any product providing higher utility than the hurdle one.
\subsection{Conjoint Analysis and Polyhedral Estimation Method}
Model SOC-D is established with a prerequisite that all $\boldsymbol{u}_n$ has already been collected. As the utilities are not directly observable metrics, we require a precedent procedure of collecting some primitive data and estimating the utilities from these data. Conjoint analysis (CA) is one of the most commonly applied approaches \citep{Wang2012A,wang2022robust}. %The methodology is composed of two parts, collecting the survey data from customers and estimating utilities from the survey data. 

\begin{figure}[htbp]
	\centering
	\includegraphics[width=0.7\linewidth]{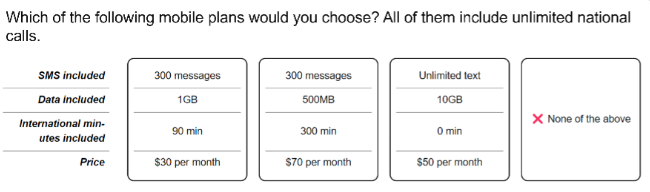}
	\caption{An example of choice-based CA on the mobile plans}
	\label{fig:cbc}
\end{figure}

In CA, the survey data are collected by asking the sampled customers to compare products with different attribute levels in the format of rating, ranking or choosing from a pre-specified set of products. Figure \ref{fig:cbc} presented an example of choice-based CA (CBC) survey question provided by a network-based CA tool, Conjointly \footnote{https://conjointly.com/solutions/feature-selection/}, where the respondents are asked to choose from different mobile plans or nothing.  The survey are usually repeated several rounds to collect enough data to estimate the partworth utilities for each sampled customer.

With the collected data, the next step is to estimate the partworth. Polyhedral method proposed by \cite{Toubia2003Fast,toubia2004polyhedral} is one of the most popular ones. In the method, the CA survey are applied to characterize a polyhedral set of utilities consistent with the response. 
For example, if the first product in Figure \ref{fig:cbc} is chosen by the customer, under the utility maximization framework, the utility consistent with the response is characterized by the linear constraints that the summation of partworth of levels selected for the first product is higher than those of the rest ones and the outside options. This can be expressed by a series of pairwise comparison between the products. Suppose the number of pairwise comparison for customer $n$ is $K_n$ and let $q_{nk}$ record the gap between the difference between the compared pairs. 
\begin{itemize}
	\item $q_{nkl}=1$ if level $l$ is selected for the first profile but not for the second one.
	\item $q_{nkl}=0$ if level $l$ is selected for both or neither profile.
	\item $q_{nkl}=-1$ if level $l$ is selected for the second profile but not for the first one.
\end{itemize}
Let $r_{nk}$ be a metric of ``intensity'' to which customer $n$ prefers the first product than the second one of comparison $k$. It can be applied to represent the ratings given by the respondents, price gaps between the chosen product and other alternatives, or a constant representing the minimal difference in utility that can be tell by the customer, etc. 

According to previous notations, an inequality $\boldsymbol{u}_n^T\boldsymbol{q}_{nk}\geq \boldsymbol{r}_n$ corresponds to the $k$th pairwise comparison for customer $n$ and all the partworth utility consistent with all $K_n$ comparisons can be expressed as follows. 
\begin{align}
	\mathcal{U}_n(\boldsymbol{Q}_n, \boldsymbol{\hat{r}}_n)\equiv\{\boldsymbol{{u}}_n|\boldsymbol{Q}_n\boldsymbol{{u}}_n\geq \boldsymbol{\hat{r}}_n\}. \label{ambiguity set}
\end{align}

\begin{figure}[htbp]
	\centering
	\includegraphics[width=0.5\linewidth]{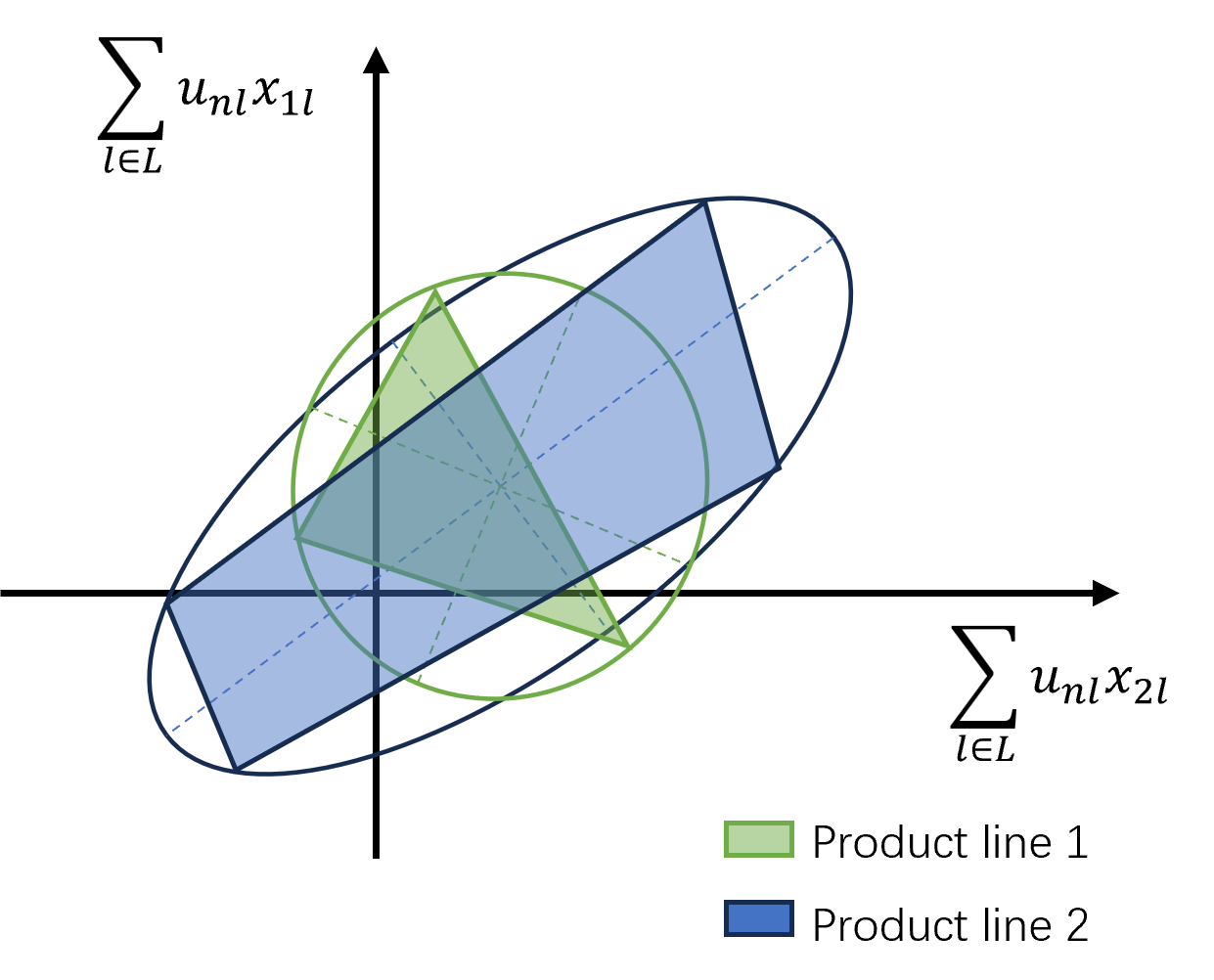}
	\caption{Drawback of neglecting the primitive survey data.}
	\label{fig:analyticcenterestimation}
\end{figure}

When the above polyhedral feasible set of $\hat{\boldsymbol{u}}$ are not shrunk to a singleton after finishing the CA survey, any point in the remaining polyhedron is consistent with the survey data. However, The classic SOC model requires a point estimation of $\boldsymbol{\hat{u}}_n$.  \cite{Toubia2003Fast} suggested to apply the analytical centers of the polyhedrons as the point estimations. We present an illustrative example of potential drawbacks when applying the analytic center estimation to the SOC model in Figure \ref{fig:analyticcenterestimation}.

Suppose we have a case with $M=2$ and two alter product lines. Figure \ref{fig:analyticcenterestimation} shows the feasible utility regions of two product lines for customer $n$ for all $\hat{u}_n\in \mathcal{U}_n(\boldsymbol{Q}_n, \boldsymbol{\hat{r}}_n)$.  The two product line are the same when evaluating against the analytical center. However, product line 1 is better with respect to the SOC problem because at least one product in product line 1 provides non-negative utility for all $\hat{u}_n\in \mathcal{U}_n(\boldsymbol{Q}_n, \boldsymbol{\hat{r}}_n)$ while there exists some  $\hat{u}_n\in \mathcal{U}_n(\boldsymbol{Q}_n, \boldsymbol{\hat{r}}_n)$ such that both products in product line 2 provide negative utility.

Because the problem does not have close-form solution, it is hard to provide the linkage between the primitive CA data to the estimated utility and ultimately to the final decision of the recommended product line, which undermines the explainability of the whole process.

In the next section, we settle the above two drawbacks of applying the analytic center estimations by using the worst-case utility over $\mathcal{U}_n(\boldsymbol{Q}_n, \boldsymbol{\hat{r}}_n)$. As we will show, the primitive data is directly taken into optimization through the dual problem of pursuing the worst-case utility, which establish the relationship between the recommended product line and the primitive CA data.

\begin{remark}[Dealing with Individual Transaction Data]
	Although collected in different setting, transaction data on individual level can be dealt in the same way as the CBC data and therefore, can be used to construct the polyhedral sets \eqref{ambiguity set}.  We refer the readers to \cite{chen2023model} for a thorough discussion in dealing with the transaction data.
\end{remark}
\begin{remark}[Response Error]
	In this paper, we focus on the case with no response error, i.e., customer always choose the alternative with highest utility. Many techniques are proposed to enhance the ability of polyhedral method to incorporate response error, such as \cite{evgeniou2005generalized,Bertsimas2013Learning}. Note that these adjustments do not influence the establishment of our model. As such, we do not discuss much on the response error issue.
\end{remark}
\section{End-to-end SOC Models}\label{sec:end-to-end}
In contrast to the classic model, we do not take the point estimation of the partworth as input. Instead, we take $\mathcal{U}_n(\boldsymbol{Q}_n,\boldsymbol{\hat{r}}_n)$\footnote{	More advanced method such as ellipsoidal method proposed by \cite{saure2019ellipsoidal} can also be applied to establish a similar model to obtain robustness against utility estimation. However, the major drawback is losing the explainability provided by the polyhedron method.} into consideration through maximizing the worst-case SOC over the polyhedrons. Formally, the proposed model is given as follows. % 
\begin{align}
	\mbox{[PCO-R]} \max_{\boldsymbol{X}\in \mathcal{X}}& \frac{1}{N} \sum_{n\in [N]} \inf_{\boldsymbol{{u}}_n \in \mathcal{U}_n(\boldsymbol{Q}_n,\boldsymbol{\hat{r}}_n)} \boldsymbol{1}\{\max_{m\in[M]} \boldsymbol{{u}}_n^T\boldsymbol{x}_m\geq h_n\} \label{PCO-R}
\end{align}
The motivation of applying the worst-case utility corresponds to the two drawbacks of applying the analytic center as point estimation. First, the approach provides robustness against uncertainty in the utility estimation. As we will show in Section \ref{sec_NE}, applying the worst-case utility outperform using the analytical center as point estimation when the number of question is small. Second, as we will show in this section, by taking the dual problem of pursuing the worst-case utility, we integrate the primitive data into the SOC optimization, which allows the model attributing the optimal product line to the primitive survey data. 

The max-min format is intractable. Based on the dual problem of pursuing the worst-case utility over $\mathcal{U}_n(\boldsymbol{Q}_n,\boldsymbol{\hat{r}}_n)$ for a given $\boldsymbol{X}$, problem (\ref{PCO-R}) is equivalent to the following tractable reformulation.
\begin{subequations} \label{problem re}
	\begin{align}
		\mbox{\rm{[PCO-RT]}} \quad & \max_{\boldsymbol{X}\in \mathcal{X}} \frac{1}{N}\sum_{n\in [N]}y_n \label{obja}\\
		\rm{s.t.} \quad & C(y_n-1)\leq \boldsymbol{\beta}_n^T\boldsymbol{\hat{r}}_n, \quad  \forall n\in [N],\label{E1}\\
		& \boldsymbol{\beta}_n^T\boldsymbol{Q}_n-\boldsymbol{\alpha}^T_n\boldsymbol{X}\leq 1-y_n,\quad  \forall n \in [N],\label{E2}\\
		& \boldsymbol{\beta}_n^T\boldsymbol{Q}_n-\boldsymbol{\alpha}^T_n\boldsymbol{X}\geq y_n-1,\quad  \forall n \in [N],\label{E2_2}\\
		& \boldsymbol{e}^T\boldsymbol{\alpha}_n= 1,\quad  \forall n\in [N],\label{E3}\\
		& \boldsymbol{\alpha}_n,\boldsymbol{\beta}_n\geq 0,\quad  \forall n\in [N],m\in [M]\label{Ee}\\
		& \sum_{l\in\mathcal{L}_a}x_{ml}=1,\forall a\in[A],m\in[M],\\
		&x_{ml},y_n \in \left\{ 0,1 \right\},\quad  \forall n\in [N], l\in[L],m\in[M].
	\end{align}
\end{subequations}
where $C$ is a large enough positive constant and $\boldsymbol{\alpha}_n,\boldsymbol{\beta}_n$ are the dual variables of pursuing the worst-case utility of customer $n$. 

Through constraint \eqref{E2}, model PCO-RT attributes the recommended product line to the primitive data. Intuitively, the attribution is achieved by recombining the survey questions to seek the evidence on the fact that the worst-case utility is non-negative. Because $\boldsymbol{Q}\boldsymbol{u}_n\geq \hat{\boldsymbol{r}}_n$, we have $\boldsymbol{\beta}_n^T\boldsymbol{Q}\boldsymbol{u}_n\geq \boldsymbol{\beta}_n^T\hat{\boldsymbol{r}}_n$. Constraint \eqref{E1} implies that $y_n=1$ only if $\boldsymbol{\beta}_n^T\boldsymbol{r}_n\geq 0$. As such, $\boldsymbol{\beta}_n^T\boldsymbol{Q}\boldsymbol{u}_n\geq 0$. For $n$ with $y_n=1$, constraint \eqref{E2_2} forces the equivalence between $\boldsymbol{\beta}_n^T\boldsymbol{Q}_n$ and $\boldsymbol{\alpha}^T_n\boldsymbol{X}$. Accordingly, $\boldsymbol{\alpha}^T_n\boldsymbol{X}\boldsymbol{u}_n=\sum_{m\in[M]}\alpha_{nm}\boldsymbol{u}_n^T\boldsymbol{x}_m\geq 0$. The left hand side of the inequality is a weighted average utility over all offered product due to the constraints on the value of $\alpha_n$. Since $\max_{m\in[M]}\boldsymbol{u}^T_n\boldsymbol{x}_m\geq \sum_{m\in[M]}\alpha_{nm}\boldsymbol{u}_n^T\boldsymbol{x}_m$, the highest utility among the $M$ products must be non-negative for all $\boldsymbol{u}_n\in\mathcal{U}_n(\boldsymbol{Q}_n,\boldsymbol{r}_n)$.

Two cases are nested when $y_n=0$. The first is when $\boldsymbol{\beta}_n^T\hat{\boldsymbol{r}}_n\leq 0$, i.e., the worst-case utility for all product negative. The other is when  $\boldsymbol{\beta}_{n}^T\hat{\boldsymbol{r}}_n$ there are no $\boldsymbol{\alpha}_n$ and $\boldsymbol{\beta}_n$ such that  $\boldsymbol{\beta}_n^T\boldsymbol{Q}_n=\boldsymbol{\alpha}^T_n\boldsymbol{X}$. This case implies that the worst-case utility is unbounded. One of the typical instances is that when an level untested in the previous questions, i.e., that is, there exist $l\in [L]$, $q_{nkl}=0,\forall n\in[N], k\in[K]$, is selected by all the products in the product line. This implies an advantage of taking the primitive data into optimization. Specifically, the model will automatically filter out the levels with no supportive evidence on its contribution to customer selection

%The rationale of involving such a constraint is twofold. First, the maximal utility for the product line is also unbounded, which implies that there are no information to infer the customer utility on the direction and therefore, establishing such a product line based on no evidence. Second, the maximal utility is bounded, which indicates that the analytic center of the utility is negative and as such, such a product line will not be selected by the nominal model. Furthermore, $|\beta_{nk}|,\forall n\in[N]$ indicates the extent of the optimal product relying on question $k,\forall k\in [K_n]$. \alert{Adding more illustration on how to use the values.}

Although the above reformulation is informative, the multiplication of the variables, i.e.,  $\boldsymbol{\alpha}_n^T\boldsymbol{X}$, impedes its tractability. We apply variables $\gamma_{nml}$ to further linearize $\alpha_{nm}x_{ml}$ to the following mixed integer linear programming (MILP) model. 
\begin{subequations} 
	\begin{align}
		\mbox{\rm{[PCO-RT]}} \quad & \max_{x_{ml},\beta_{nk},\alpha_{nm}} \frac{1}{N}\sum_{n\in [N]}y_n \\
		\rm{s.t.} \quad & C(y_n-1)\leq \sum_{k}\beta_{nk}\hat{r}_{nk}, \quad  \forall n\in [N],\label{purchase}\\
		&\sum_{k\in[K]}\beta_{nk}q_{kl}-\sum_{m\in[M]}\gamma_{nml}\leq 1-y_n,\quad  \forall n \in [N],l\in[L],\label{ub_linearcombination}\\
		&\sum_{k\in[K]}\beta_{nk}q_{kl}-\sum_{m\in[M]}\gamma_{nml}\geq y_n-1,\quad  \forall n \in [N],l\in[L],\\
		&\gamma_{nml}\leq x_{ml},\forall n\in [N],m\in [M],l\in [L],\label{gamma1}\\
		&\sum_{l\in \mathbb{L}_a}\gamma_{nml}= \alpha_{nm},\forall n\in [N],m\in [M],a\in [A],\label{gamma2}\\
		&\sum_{m\in[M]}\alpha_{nm}= 1,\quad  \forall n\in [N]\\
		&\gamma_{nml},\alpha_{nm},\beta_n\geq 0,\quad  \forall n\in [N],m\in [M],l\in[L]\\
		&y_n \in \left\{ 0,1 \right\},\quad  \forall n\in [N],
	\end{align}
\end{subequations}
Constraint \eqref{gamma1} and \eqref{gamma2} enforce the equivalence between $\gamma_{nml}$ and $\alpha_{nm}$. Constraints \eqref{gamma1} forces $\gamma_{nml}=0$ when $x_{ml}=0$, which is equal to $\alpha_{nm}x_{ml}$. For $x_{ml}=1$, we need to force $\gamma_{nml}=\alpha_{nm}$. Constraint \eqref{gamma2} force the equivalence by taking advantages of the constraint $\sum_{l\in \mathcal{L}_a}x_{ml}=1$. Since $x_{ml}=1$ for one and only one level over $\mathcal{L}_a$, exactly one $\gamma_{nml}$ is allowed to be non-zero over all $l\in\mathcal{L}_a$. As such, constraint \eqref{gamma2} enforces $\gamma_{nml}=\alpha_{nm}$ for $x_{ml}=1$.

\subsection{Performance Guarantee}
Except for its explainability, applying the worst-case utility allows us to characterize the performance of the model as a function of both estimation accuracy and sample size. %We derive the asymptotic performance guarantee of the proposed model with bounded polyhedron and the following mild assumption on the Lipchitz continuity of the cumulative density function (CDF) of the random variable $\max_{m\in[M]} {\boldsymbol{\tilde{u}}}^T{\boldsymbol{x}_m},\forall \boldsymbol{X}\in\mathcal{X}$, in which $\boldsymbol{\tilde{u}}$ follows the ground-truth distribution, denoted as $\Pi^*$, generating the oracle utility of the sampled customers.

Let $\boldsymbol{X}^N$ be the optimal solution to the proposed model with rows  and $\boldsymbol{X}^*$ be the oracle optimal solution defined as
\begin{eqnarray}
	\boldsymbol{X}^*\equiv\arg\max_{\boldsymbol{X}\in\mathcal{X}} E_{\Pi^*}\{1(\max_{m\in[M]}\boldsymbol{\tilde{u}}^T\boldsymbol{x}\geq 0)\}
\end{eqnarray}
where $\boldsymbol{\tilde{u}}$ follows the ground-truth distribution, denoted as $\Pi^*$, generating the oracle utility of the sampled customers. Let $\hat{u}_n,\forall n\in[N]$ be the unobserved ground-truth partworth of the sample customer $n$.

%First of all, we present the asymptotic performance guarantee of the SOC model with the known oracle utilities of the sampled customers as a benchmark.  
%\begin{proposition}\label{pro:outofsample2}
%	Suppose $\{\boldsymbol{\hat{u}}_n\}_{n=1}^N$ is i.i.d with a common distribution $\boldsymbol{\tilde{u}}$. 
%	Then for arbitrary $\epsilon>0$ and $\delta\in (0,1)$.
%	We have that with probability higher than $1-\delta$,
%	\begin{align*}
%		\left|P(\max_{m\in[M]}\boldsymbol{\tilde{u}}^T\boldsymbol{x}^{N}_m)-P(\max_{m\in[M]}\boldsymbol{\tilde{u}}^T\boldsymbol{x}_m^*)\right|\leq\epsilon
%	\end{align*}
%	when $N \geq \frac{2e^2(ML+2)\textnormal{log}2}{\epsilon^2}\textnormal{log}(\frac{2}{\delta})$ where $e$ is the Euler's number.
%\end{proposition}

In Proposition \ref{pro:outofsample1}, we present the performance guarantee of the propose model as a function of the maximal size of polyhedrons, the dimension of the attributes, the number of products, and the sample size.
\begin{proposition}\label{pro:outofsample1}
	Suppose that 
	\begin{itemize}
		\item [a)] the distance between any two points in $\mathcal{U}_n(\boldsymbol{Q}_n,\boldsymbol{\hat{r}}_n),\forall n\in[N]$ is no more than $d,\forall n\in[N]$
		\item [b)] and for every $\boldsymbol{X}$, the CDF of $ \max_{m\in[M]}\tilde{u}^T \boldsymbol{x}_m$ is Lipschitz continuous with Lipschitz constant $\theta$, i.e., $| P(\max_{m\in[M]} \boldsymbol{u}_1^T\boldsymbol{x}_m)-P(\max_{m\in[M]}\boldsymbol{u}_2^T\boldsymbol{x}_m\geq 0)|\leq \theta|\boldsymbol{u}_1^T-\boldsymbol{u}_2^T|$
	\end{itemize}
	Then, 
	\begin{align*}
		\mathbb{E}\left| P( \max_{m\in[M]} {\boldsymbol{\tilde{u}}}^T{{\boldsymbol{x}}_m^N}\geq 0)-P( \max_{m\in[M]} {\boldsymbol{\tilde{u}}}^T{{\boldsymbol{x}}_m^*}\geq 0) \right|\leq 4{\textnormal{log}(2)} e\sqrt{2(LM+2)/N}+\theta \sqrt{A} d,
	\end{align*}
	where $e$ is the Euler's number.
\end{proposition}
The performance guarantee is composed of two parts. 
\begin{itemize}
	\item The first part is a common term when SAA is applied to approximate the oracle problem. The term reflects the uncertainty caused by applying the sampled distribution to the population one. As sample size increases, this term diminishes in a rate of square-root. We can also analyse the tail probability of the absolute difference between the data-driven and true optimal market-share, which we summarize in appendix.
	\item The second part is related to the overall estimation accuracy among the sample customers. When the CA questions are properly designed, $d$ reduces in the number of questions is common in previous literature. For example, \cite{toubia2004polyhedral} started with a bounded polyhedral set that all the partworths are lower bounded and their summation are upper bounded and propose an ACA framework to halved the longest distance between any two points in each iteration. With this ACA framework, the original $d$ is cut into half and therefore $d$ will keep reducing after each round. Accordingly, the bound in Proposition \ref{pro:outofsample1} get tighteras the number of question increases.
\end{itemize}

\section{Adaptive Conjoint Analysis based on Column Generation}\label{sec:ACA}
Design the CA survey adaptively based on the responses to previous questions is a major focus of the CA studies. This scheme is known as ACA for short. In most of the previous literature \citep{Toubia2003Fast,toubia2004polyhedral,saure2019ellipsoidal}, the ACA approaches are proposed to reduce the estimation uncertainty as soon as possible. For example, \cite{Toubia2003Fast,toubia2004polyhedral} proposed to ask a new question corresponding to the vector most nearly parallel to the longest axis of the circumscribed ellipsoid of the polyhedral set. However, the information contained in a question most efficiently reducing estimation uncertainty is not necessary the one that benefits most in the aspect of improving the objective of the PLD optimization \citep{JooOptimal}. 

As such, leveraging on the structure,we propose a goal-directed ACA approach to align the CA design with the optimization of model PCO-RT. Except for its explainability, model PCO-RT indicates how the asked questions and response contributing to the optimization. Since the asked questions corresponds to a series of additional variables $\beta_{nk},\forall n\in[N],k\in[K]$, asking a new question  $n$, $q_{n,K_n+1}$ to sampled customer $n$ adds a new variable $\beta_{n,K+1}$ to the programming model. This is with the same idea of column generation, which start with a restricted master problem (RMP) that includes only a subset of all possible variables and then iteratively add new ones to improve the objective of the RMP as fast as possible. In our setting, model PCO-RT can be seen as such a RMP with $\beta_{nk}$ as the variables needed to be added to it. 

In column generation approach, a pricing sub-problem (PP) is required to find the variables that can potentially improve the objective function. Let $\lambda^*_{n},\boldsymbol{\eta}^*_{n},\boldsymbol{\kappa}^*_{n}$ denote the optimal dual variables of constraints \eqref{purchase}-\eqref{gamma1} in the linear relaxation of model PCO-RT with $\boldsymbol{X}$ fixed to the optimal solutions with the current $K$ questions. When the response to the new question, i.e., $r_{n,K+1},\forall n\in[N],$ is known, PP problem is formulated as follows
\begin{eqnarray}
	\mbox{(PP)}&\max&  \sum_{n\in[N]}\hat{r}_{n,K+1}\lambda_{n}^*+ (\boldsymbol{\eta}^*_{n}-\boldsymbol{\kappa}^*_{n})^T\boldsymbol{q}_{n,K+1}\\
	&s.t.&  \boldsymbol{q}_{n,K+1}=\boldsymbol{x}_{n1}-\boldsymbol{x}_{n2},\forall n\in [N],\\
	&& \sum_{l\in\mathcal{L}_a}x_{nml}=1,\forall a\in[A],n\in[N],m=\{1,2\}
\end{eqnarray}
where $\boldsymbol{x}_{n1},\boldsymbol{x}_{n2}$ are the alternatives compared in the $K+1$ round.
However, the above model is not obtainable because we can not know the exact response of the customers before the question is asked. For example, if the question is generated by a metric CA, i.e., $\hat{r}_{n,K+1}$ is exactly the utility gap between the two compared products, $\hat{r}_{nk}$ can be any value in the projected polyhedron $\{\sum_{l\in[L]}u_{l}q_{kl}|\boldsymbol{u}\in \mathcal{U}(\boldsymbol{Q}_n,\hat{\boldsymbol{r}}_n)\}$; if the question is generated by choice-based CA between two alternatives $\boldsymbol{x}_1,\boldsymbol{x}_2$,  $\boldsymbol{q}_{k}$ can be either $\boldsymbol{x}_1-\boldsymbol{x}_2$ or $\boldsymbol{x}_2-\boldsymbol{x}_1$ according to the customer choice. 

To tackle the issue, we design a novel approach to predict the response by utilizing the CA survey data on an aggregated level. Specifically, we estimate on the population-level partworth utility distribution and then use the conditional probability or expectation on the observed data to establish model PP for each respondent. We illustrate the approach for both metric CA and CBC in the following two sections, respectively. In this section, we assume that the sampled utility are i.i.d samples from a multivariate normal distribution.
\subsection{Metric CA}
For metric CA, the responses are exactly the gaps between the compared products and therefore, equations $\boldsymbol{Q}_n\boldsymbol{u}_n=\hat{\boldsymbol{r}}_n$ are applied to substitute the inequalities \citep{Toubia2003Fast}. Under this assumption, $r_{n,K+1}=\boldsymbol{u}_n^T\boldsymbol{q}_{n,K+1}$. Accordingly, model PP is presented as follows in this situation.
\begin{eqnarray*}
	\mbox{(PP-M)}&\max&  \sum_{n\in[N]} (\lambda_{n}^*\boldsymbol{u}^T_n+\boldsymbol{\eta}^*_{n}-\boldsymbol{\kappa}^*_{n})\boldsymbol{q}_{n,K+1}\\
	&s.t.&  \boldsymbol{q}_{n,K+1}=\boldsymbol{x}_{n1}-\boldsymbol{x}_{n2},\forall n\in [N],\\
	&& \sum_{l\in\mathcal{L}_a}x_{nml}=1,\forall a\in[A],n\in[N],m=\{1,2\}
\end{eqnarray*}
Still, the main difficulty of establishing the above model is the uncertainty in $\boldsymbol{u}_n$. Suppose that $u_n$ are i.i.d sample from a multivariate normal distribution controlled by mean vector $\boldsymbol{\mu}$ and covariance matrix $\boldsymbol{\Sigma}$. We can first use the response of all customers to estimate the parameters of the oracle distribution and then apply the conditional distribution on the observed data $(\boldsymbol{Q}_n,\boldsymbol{r}_n),\forall n\in[N]$, $P_{\boldsymbol{u}|\boldsymbol{Q}^T_nu=\boldsymbol{r}_n}$, to get more information on the value of $\boldsymbol{{u}}_n$. With the conditional probability, we maximize the expected reduced cost as follows to get the next question.
\begin{eqnarray*}
	\mbox{(PP-EM)}&\max&  \sum_{n\in[N]}E_{P_{\boldsymbol{u}|\boldsymbol{Q}^T_nu=\boldsymbol{r}_n}}(\sum_{l\in[L]}(u_{l}\lambda_{n}^*+ \eta^*_{nl}-\kappa^*_{nl}) q_{K+1,l})\\
	&s.t.&  \boldsymbol{q}_{n,K+1}=\boldsymbol{x}_{n1}-\boldsymbol{x}_{n2},\forall n\in [N],\\
	&& \sum_{l\in\mathcal{L}_a}x_{nml}=1,\forall a\in[A],n\in[N],m=\{1,2\}
\end{eqnarray*}

%Many parametric and non-parametric approaches can be applied to estimate the distributions and optimize model PP-EM. When $u_n$ are i.i.d sample from a multivariate normal distribution controlled by mean vector $\boldsymbol{\mu}$ and covariance matrix $\boldsymbol{\Sigma}$, we can solve the problem in a more neat way. 

Denote the estimated mean and covariance matrix as $\hat{\boldsymbol{\mu}}_N$ and $\hat{\boldsymbol{\Sigma}}_N$, we can further reformulate PP-EM as follows.
\begin{eqnarray*}
	\mbox{(PP-EMT)}&\min&  \sum_{n\in[N]} \lambda_n^*(\boldsymbol{q}_{n,K+1}^T\hat{\boldsymbol{\mu}}_N+\boldsymbol{q}_{n,K+1,}^T\hat{\boldsymbol{\Sigma}}_N\boldsymbol{Q}^T_n(\boldsymbol{Q}_n\hat{\boldsymbol{\Sigma}}_N\boldsymbol{Q}_n^T)^{-1}(\boldsymbol{r}_n-\boldsymbol{Q}_n\hat{\boldsymbol{\mu}}_N))+\boldsymbol{q}_{n,K+1}^T(\boldsymbol{\eta}^*_n-\boldsymbol{\kappa}^*_n)\\
	&s.t.&  \boldsymbol{q}_{n,K+1}=\boldsymbol{x}_{n1}-\boldsymbol{x}_{n2},\forall n\in [N],\\
	&& \sum_{l\in\mathcal{L}_a}x_{nml}=1,\forall a\in[A],n\in[N],m=\{1,2\}
\end{eqnarray*}
Although the above model seems quite complex, the above model is a mixed integer linear programming because only $\boldsymbol{q}_{n,K+1},\boldsymbol{x}_{n1},\boldsymbol{x}_{n2}$ are decision variables. %Furthermore, because one equality constraint will reduce one dimension of the polyhedron $\mathcal{U}_n(\boldsymbol{Q}_n,\boldsymbol{\hat{r}}_n)$, we let $\mbox{rank}(\boldsymbol{Q}_n)<L$. Otherwise,  $\mathcal{U}_n(\boldsymbol{Q}_n,\boldsymbol{\hat{r}}_n)$ will be a singleton. If the rows of $\boldsymbol{Q}_n$ are linearly independent, $\mbox{rank}(\boldsymbol{Q}_n)=K$ if $K\leq L$. In this situation, $\boldsymbol{Q}\boldsymbol{Q}^T_n$ is invertable and therefore, $\boldsymbol{Q}_n\boldsymbol{\Sigma}\boldsymbol{Q}_n^T$ is invertable, which 

%Secondly, the estimation can be dealt in a simple way. Since each $\boldsymbol{r}_n=\boldsymbol{Q}_n\boldsymbol{u}$, $\boldsymbol{r}_n$, stack all the observation data to form the following equation:
%\[ \mathbf{R} = \mathbf{Q} \mathbf{u}\]
%where \(\mathbf{R} = [\mathbf{r}_1, \mathbf{r}_2, \ldots, \mathbf{r}_n]^T\) and \(\mathbf{Q} = [\mathbf{Q}_1, \mathbf{Q}_2, \ldots, \mathbf{Q}_n]^T\). And then, we can use the pseudo-inverse to solve the least squares estimation of the realization of the partworth utilities.
%\begin{equation}
%	\hat{\mathbf{u}} = (\mathbf{Q}^T \mathbf{Q})^{-1} \mathbf{Q}^T \mathbf{R}\label{eq:inferedhatu}
%\end{equation}
%Consequently, we can estimate the mean vector and covariance matrix of \(\mathbf{u}\) by the sampled mean $\hat{\boldsymbol{\mu}}_N$ and covariance $\hat{\boldsymbol{\Sigma}}_N$. 
\subsection{Choice-based CA}
Choice-based CA is different to the metric one in two aspects. First, the value of $r_n$ is not the exact utility gap between the compared alternatives. It is either be 0 or determined by the gap of some metric attributes such as price. Secondly, the choice determines the specific value of the questions. For example, when offering two products, says $\boldsymbol{x}_1$ and $\boldsymbol{x}_2$, the question can be either $\boldsymbol{x}_1-\boldsymbol{x}_2$ if $\boldsymbol{x}_1$ is chosen or  $\boldsymbol{x}_2-\boldsymbol{x}_1$ otherwise.  

In this section, we focus on a case where customer $n$ choose from a pair of products, $\boldsymbol{x}_{n1}$ and $\boldsymbol{x}_{n2}$ and the response $r_{n,K+1}$ measured by the gap between the prices associated with the products, denote as $p_{1k}$ and $p_{2k}$, respectively.  Then, with the conditional expected reduced cost as objective, the new question can be found by the following problem.
\begin{footnotesize}
	\begin{eqnarray*}
		\mbox{(PP-C)}&\max&  \sum_{n\in[N]}P\left(\sum_{l\in[L]}u_l(x_{l1}-x_{l2})-(p_1-p_2)\geq 0|\boldsymbol{Q}^T_nu\geq \boldsymbol{r}_n\right)\left((p_1-p_2)\lambda_{n}^*+\sum_{l\in[L]}( \eta^*_{nl}-\kappa^*_{nl})(x_{l1}-x_{l2})\right)\\
		&&+P\left(\sum_{l\in[L]}u_l(x_{l1}-x_{l2})-(p_1-p_2)\leq 0|\boldsymbol{Q}^T_nu\geq \boldsymbol{r}_n\right)\left((p_2-p_1)\lambda_{n}^*+\sum_{l\in[L]}( \eta^*_{nl}-\kappa^*_{nl})(x_{l2}-x_{l1})\right)\\
		&s.t.& \sum_{l\in\mathcal{L}_a}x_{nml}=1,\forall a\in[A],n\in[N],m=\{1,2\}\\
		&& p_1,p_2\in\mathcal{P},x_{ml}\in\{0,1\}
	\end{eqnarray*}
\end{footnotesize}
The calculation of the conditional expected improvement is composed of two parts. Specifically, the first and second items of the objective represents the situation where the first and second option are chosen, respectively.

The above problem is more difficult compared to that for the metric CA from the computational tractability. Even if the utilities are i.i.d samples from the multivariate normal distribution, there are no closed-form for the conditional expectation. To deal with the issue, we use an SAA approach by draw $S$ samples from the estimated conditional distribution and apply the frequency over the samples as a proxy of the probability. Denote the $s$th sample from the conditional distribution by $\boldsymbol{\hat{u}}_{ns}\sim P\{\boldsymbol{u}\sim N(\hat{\boldsymbol{\mu}},\hat{\boldsymbol{\Sigma}})|\boldsymbol{Q}^T_nu=\boldsymbol{r}_n\}$. Then, we approximate model PP-C by the following model.

\begin{footnotesize}
\begin{eqnarray*}
	\mbox{(PP-CT)}&\max&  \sum_{n\in[N]}\frac{1}{S}\sum_{s\in[S]} y_{ns}\left((p_1-p_2)\lambda_{n}^*+\sum_{l\in[L]}( \eta^*_{nl}-\kappa^*_{nl})(x_{l1}-x_{l2})\right)+(1-y_{ns}) \left((p_2-p_1)\lambda_{n}^*+\sum_{l\in[L]}( \eta^*_{nl}-\kappa^*_{nl})(x_{l2}-x_{l1})\right)\\
	&s.t.&  \sum_{l\in\mathcal{L}_a}x_{ml}=1,\forall a\in[A],m=\{1,2\}\\
	&&\hat{\boldsymbol{u}}_{ns}^T(\boldsymbol{x}_1-\boldsymbol{x}_2)-(p_1-p_2)\geq C(y_{ns}-1)\\
	&& p_1,p_2\in\mathcal{P},x_{ml}\in\{0,1\}
\end{eqnarray*}
\end{footnotesize}
A binary variable $y_{ns}$ is applied to indicate whether the sample $s$ for calculating the conditional probability of customer $n$ satisfying $\hat{\boldsymbol{u}}_{ns}^T(\boldsymbol{x}_1-\boldsymbol{x}_2)-(p_1-p_2)\geq 0$. Further linearizing and simplifying the model, we have 

\begin{eqnarray*}
	\mbox{(PP-CTL)}&\max_{x,p,y,z}& \sum_{n\in[N]}\sum_{s\in[S]}z_{ns}\\
	&s.t.&  \sum_{l\in\mathcal{L}_a}x_{ml}=1,\forall a\in[A],m=\{1,2\}\\
	&&\sum_{l\in[L]}\hat{u}_{nsl}(x_{1l}-x_{2l})-(p_1-p_2)\geq C(y^+_{ns}-1)\\
	&&y_{ns}^++y_{ns}^-=1,\forall n\in[N],s\in[S],\\
	&&z_{ns} \leq \left((p_1-p_2)\lambda_{n}^*+\sum_{l\in[L]}( \eta^*_{nl}-\kappa^*_{nl})(x_{l1}-x_{l2})\right)+C(1-y^+_{ns})\\
	&&z_{ns} \leq -\left((p_1-p_2)\lambda_{n}^*+\sum_{l\in[L]}( \eta^*_{nl}-\kappa^*_{nl})(x_{l1}-x_{l2})\right)+C(1-y^-_{ns})\\
	&& p_1,p_2\in[\underline{p},\bar{p}),x_{ml}\in\{0,1\}
\end{eqnarray*}

%Another difficulty is the estimation of partworth distribution on an aggregated level since we can not apply the pseudo-inverse to do the estimation. Instead, in the case where the same comparison pair $\boldsymbol{x}_{1k},\boldsymbol{x}_{2k}$ and prices $p_{1k},p_{2k}$ are asked to each customer in each question $k,\forall k\in[K]$, we construct the least square estimation by the share of customer choosing $\boldsymbol{x}_{1k}$ over $N$ customers. We denote this share by $\hat{P}_N((\boldsymbol{x}_{1k}-\boldsymbol{x}_{1k})^T\boldsymbol{u}\geq p_{1k}-p_{2k})$ and solving the following problem.
%\begin{eqnarray}
%	\max_{\hat{\mu},\hat{\boldsymbol{\Sigma}}}\sum_{k\in[K]}\left(\hat{P}_N\left((\boldsymbol{x}_{1k}-\boldsymbol{x}_{1k})^T\boldsymbol{u}\geq p_{1k}-p_{2k}\right)-\Phi\left(\frac{p_{1k}-p_{2k}-(\boldsymbol{x}_{1k}-\boldsymbol{x}_{1k})^T\boldsymbol{\mu}}{(\boldsymbol{x}_{1k}-\boldsymbol{x}_{1k})^T\boldsymbol{\Sigma}(\boldsymbol{x}_{1k}-\boldsymbol{x}_{1k})}\right)\right)^2
%\end{eqnarray}
%where $\Phi$ is the cdf of the standard normal distribution.
%It would be hard to directly solve the above solution due to the nonlinear conditional probability $P\left(\sum_{l\in[L]}u_l(x_{l1}-x_{l2})-(p_1-p_2)\geq 0|\boldsymbol{Q}^T_nu\geq \boldsymbol{r}_n\right)$. We apply a bootstrapping approach by sampling from the normal distribution  $\boldsymbol{u}\sim \hat{\mu},\hat{\boldsymbol{\Sigma}}$ truncated by $\boldsymbol{Q}_n\boldsymbol{u}\geq \boldsymbol{r}_n$. 

\begin{remark}
	There are two key directions where we can easily extend our ACA framework.
	\begin{itemize}
		\item Note that although the normal distribution is assumed for the ground-truth partworth utilities, it is mainly for the ease of illustration. Other assumption on the utility distribution at the population level can also be applied. 
		\item Furthermore, the ACA framework is also applicable to both cases where different customers are asked to compare the same or different pairs of problems. Both models PP-EMT and PP-CTL keep their linearity in both cases 
	\end{itemize} 
\end{remark}

\section{Numerical Experiments} \label{sec_NE}
In this section, we implement a simulation study to validate the efficiency of the proposed product line design and ACA frameworks. We benchmark our model with a typical used estimate-then-optimize framework where the point estimation of partworths are first obtained and the product line is then selected by model SOC-D with the estimation as the input. Analytical center (AC) estimation \citep{Toubia2003Fast} is applied to get the point estimations. We do not use hierarchical bayesian method since its assumes on the multinomial logit choice behavior does not aligned with the first-choice assumption behind the SOC problem.  For the survey design, a static orthogonal design generated through SPSS 25.0, whose detail is shown in Appendix \ref{OrthCA}, and the fast-polyhedral ACA (FPACA) proposed in \cite{Toubia2003Fast,toubia2004polyhedral} are considered for comparison. Consequently, four essential combinations are implemented to show the effectiveness of the proposed model and ACA framework. 
\begin{itemize}
	\item The first two settings are with the fixed orthogonal design. Model SOC-D with AC estimation and model PCO-RT are applied to optimize the product line for the first and second settings, respectively. The comparison of these two will shows the relative performance of the proposed model and the estimate-then-predict framework.
	\item The third setting runs with the proposed model PCO-RT and the column-generation-based ACA framework. Comparing it with the second setting isolates the benefits of the proposed ACA framework.
	\item The fourth setting implements model SOC-D with FPACA design. The comparison between the third and fourth setting demonstrates the overall benefits obtained from the proposed framework.
\end{itemize} 

All experiments are coded with python and solved Gurobi with default settings and run on a computer with a CentOS 7.9 operation system, 128 GB RAM, and Intel (R) Xeon(R) X5660 @2.80 GHz CPU.
\subsection{Data Generation Scheme}
We consider a synthetic case where there are 31 levels grouped into 4 attributes. The numbers of levels categorized into the four attributes are $5,8,9,9$, respectively. For the ease of solving, we suppose the data are collected through metric CA. The true partworth of the respondents for each instance is generated in a two-step procedure.
\begin{itemize}
	\item Firstly, we draw the mean of the partworth for each level uniformly between 4 and -5 and standard deviation equal to 10, which guarantee the customer preference disperse enough and therefore, there are no levels strictly preferred by the whole population. Further, through this randomization, we avoid the case where the orthogonal design happens to perform good or bad by coincidence. 100 instances of the parameters are established.
	\item Then, for each mean and variance sampled in the first step, we draw the individual  partworths of respondents from a multivariate normal distribution with the mean and variance. These ground-truth partworth are applied to generate the response to the CA questions and calculate the true SOC within the sampled respondents.
\end{itemize}  
Further, we apply a product composed of levels 1,6,14,23 as the outside option and without loss of generality, we calibrate its utility to 0. Accordingly, the partworths of 27 levels need to be estimated. Accordingly, the uncertainty set will converge to a singleton with 27 metric pairwise comparison questions if the CA survey is non-degenerate. To initialization the ACA methods, we take the first 5 questions from the orthogonal design . 

\subsection{Result}\label{EPCO}
\begin{figure}[htbp]
	
	\subfloat[Mean]{
		\begin{minipage}{0.5\textwidth}
			\includegraphics[width=\linewidth]{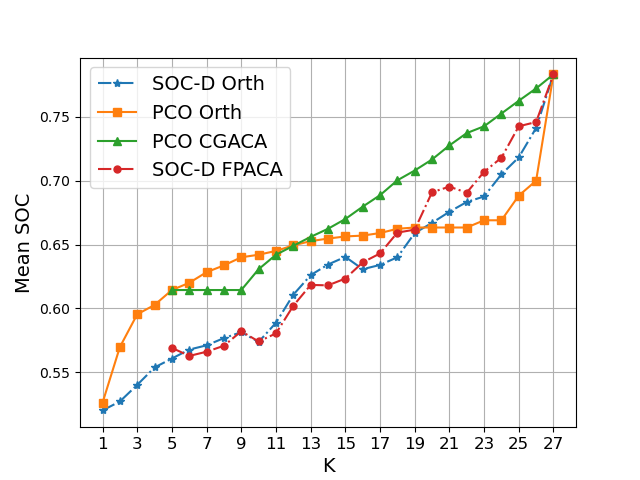}
			\label{fig:socis20}
		\end{minipage}
	}
	\subfloat[Standard Deviation]{
		\begin{minipage}{0.5\textwidth}
			\includegraphics[width=\linewidth]{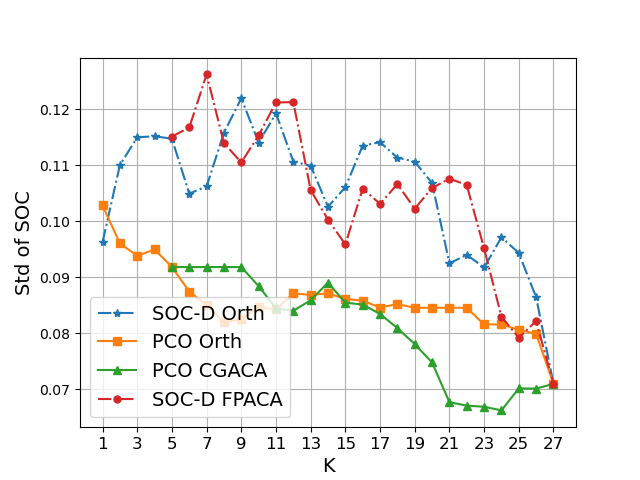}
			\label{fig:socos20}
		\end{minipage}
	}
	\caption{Comparison of average performance over model PCO-RT and model SOC-D with different estimation methods}
	\label{ModelComparison}
\end{figure}

We test the problems with different numbers of respondents and products. As all the results are similar, we exemplify the results by $N=30$ and $M=1$. For all four settings, we solve the models for the number of questions increased from 1 (start from 5 for the ACA methods) to 27 to get the corresponding product lines. Figures \ref{fig:socis20} and \ref{fig:socos20} show the mean and standard deviation of the SOC of the obtained product lines over 100 instances, respectively. 

The comparison between the first (blue lines) and second (orange lines) settings shows the outperformance of model PCO-RT to the classic estimate-then-optimize framework. As shown in Figure \ref{fig:socis20}, when the number of questions is small, the proposed model significantly improve the average performance. For example, the SOC of the product obtained by the proposed model is 11.9\% on average (from 0.573 to 0.642) when 10 questions are asked.  The observation confirms the ability of our model to achieve desirable performance when the CA survey data is not able to accurately estimate the customers' partworth. Moreover, Figure \ref{fig:socos20} indicates that the performance of the proposed model is less varied compared to that of the estimate-then-optimize, which further validates the robustness of the proposed model against partworth utility estimation uncertainty.

Comparing the fourth setting (green line) to others shows the joint benefits collected by the proposed PCO model and the column-generation-based ACA (CGACA) method. As shown in Figure \ref{fig:socis20}, with the orthogonal design, although model PCO performs outperform the classic estimate-then-optimize approach when the number of questions, the advantage of the model do not keep as the number of questions grows. In contrast, when combined with the CGACA, model PCO maintains the advantages over the estimate-then-optimize approach, which shows the importance of aligning the question design with the subsequent optimization. For example, when 21 questions are asked, model PCO with CGACA design performs 14.3\% better than model SOC-D with orthogonal design while the model with the orthogonal design performs 1.7\% worse than the SOC-D model. As shown in \cite{Toubia2003Fast}, FPACA can reduce the partworth uncertainty as fast as possible. However, as shown by Figures \ref{fig:socis20} and \ref{fig:socos20}, it has limited power in improving the average $SOC$ and reducing the variation. The comparisons among the performance of the settings implies the importance of aligning the CA design and the subsequent PLD optimization.

\section{Conclusions}\label{sec_Con}
In this paper, we propose end-to-end PLD optimization framework, where primitive CA data is directly input into the optimization Comparing to the classic ``estimate-then-optimize'' models, the end-to-end model leverages desirable explainability of how the decision relying on the survey data and allows us to present a novel performance guarantee of the model as a function of both sample size and estimation accuracy. To the best of our knowledge, the influence of the estimation accuracy on the PLD optimization is not explored before. Moreover, the structure of the model inspire us to design a goal-directed ACA framework based on the idea of column generation, which creatively takes advantages of the data on both individual and aggregated level. Through extensive numerical experiments, we show that the proposed model and the ACA framework significantly improves the performance compared to the classic ``estimate-then-optimize'' models, especially in the case with high partworth estimation uncertainty.

Future research can be extended in three directions. Firstly, except polyhedral method, how to integrate other estimation approaches, such as the hierarchal Bayesian framework \citep{allenby1998marketing}, ellipsoidal method \citep{saure2019ellipsoidal}, etc., to the subsequent optimization, such as pricing, assortment planning, and product line design, is a fruitful direction. Secondly, active treatments design in the retail context, such as recommendation and promoting, is another possible extension. As e-commerce develops, a huge amount of transaction and shopping data is available to researchers. Establishing data-driven model based on these data to find better operation strategies and characterize the customer behavior would be benefit to these companies. Third, this paper shows the possibility of enhancing explainability through robust optimization techniques. It would be interesting to investigate whether the idea can be extended to other problems.  
\bibliographystyle{pomsref} 

\let\oldbibliography\thebibliography
\renewcommand{\thebibliography}[1]{%
	\oldbibliography{#1}%
	\baselineskip14pt %Change this for line spacing within the same reference
	\setlength{\itemsep}{10pt}% %Change this for spacing between two referneces
}
\bibliography{soc}
%% Here starts the e-companion (EC)
%%%%%%%%%%%%%%%%%%%%%%%%%%%%%%%%%%%%%%%%%%%%%%%%%%%%%%%%%%
\ECSwitch

%\ECDisclaimer
%%%%%%%%%%%%%%%%%%%%%%%%%%%%%%%%%%%%%%%%%%%%%%%%%%%%%%%%%%

%%% Main head for the e-companion
\ECHead{E-companion}
\section{Proof}

\subsection{Proof of the equivalence of the reformulation.}
\textbf{Proof.} Given $\boldsymbol{X}$, the worst case utility over $\mathcal{U}_n,\forall n\in[N]$ is calculated as follows.
\begin{subequations}
	\begin{align}
		\mbox{\rm{[W-Utility]}} \quad  & \min_{\boldsymbol{u}_{n}\in\mathcal{U}_n} c_n\\
		\rm{s.t.} \quad &\sum_{l\in [L]} u_{nl}x_{ml}- c_n \leq 0, \quad \forall m\in [M]\label{W1}\\
		& \sum_{l\in L}q_{nlk}u_{nl}\geq \hat{r}_{nk}, \quad \forall k\in [K_n],\label{W2}
	\end{align}
\end{subequations}

Denote $\bm{\alpha}$ and $\bm{\beta}$ as the dual variable for the constraints (\ref{W1}) and (\ref{W2}). Take the dual of problem W-Utility, we have,
\begin{subequations}
	\begin{align}
		\mbox{\rm{[W-Utility-D]}} \quad &\max\sum_{k\in [K_n]} \hat{r}_{nk}\beta_{nk} \\
		\rm{s.t.} \quad &\sum_{k\in [K_n]}q_{nlk}\beta_{nk}-\sum_{m\in [M]} x_{ml}\alpha_{nm}=0, \quad \forall l\in [L],\label{W3}\\
		&\sum_{m\in [M]}\alpha_{nm}= 1.\label{W4}\\
		&\alpha_{nm},\beta_{nk}\geq 0, \quad \forall m\in [M],k\in[K_n]
	\end{align}
\end{subequations}

There are two cases for the problems. Firstly, model W-Utility is bounded and therefore, the two problems shares the same objective value. In this case, $\inf\limits_{\hat{u}_{nl}\in \mathcal{U}_n(\boldsymbol{Q}_n,\boldsymbol{\hat{r}}_n)}\boldsymbol{1}\{\max_{m\in[M]}\hat{u}_n^T\boldsymbol{x}_m\geq0\}=1$ if the minimal $c_n$ (also the maximal $\sum_{k\in [K_n]} \hat{r}_{nk}\beta_{nk}$) is positive and 0 otherwise. We apply a binary variable $y_n$ to replace the indicator function, whose value is determined by the following inequality, 
$$C(y_n-1)\leq \min c_n$$ 
where the minimal is taken over constraints \eqref{W1} and \eqref{W2}.
Due to the strong duality, we have
\begin{equation}
	C(y_n-1)\leq \max_{\beta_{nk}} \sum_{k\in [K_n]} \hat{r}_{nk}\beta_{nk} \label{W5}
\end{equation} 
where the feasible region of $\beta_{nk}$ is subject to constraints (\ref{W3}) and (\ref{W4}).

%Obviously, if the minimal $c_n$ is positive, $\inf\limits_{\hat{u}_{nl}\in \mathcal{U}_n(\boldsymbol{Q}_n,\boldsymbol{\hat{r}}_n)}\boldsymbol{1}\{\max_{m\in[M]}\hat{u}_n^T\boldsymbol{x}_m\geq0\}=1$. We apply a binary variable $y_n$ to replace the indicator function, whose value is determined by the following inequality, 

Inequality (\ref{W5}) has the same effect with the following set of constraints.

\begin{eqnarray}
	&&C(y_n-1)\leq \sum_{k\in [K_n]} \hat{r}_{nk}\beta_{nk} \label{W6},\\
	&&\sum_{k\in [K_n]}q_{nlk}\beta_{nk}-\sum_{m\in [M]} x_{ml}\alpha_{nm}=0, \quad \forall l\in [L],\\
	&&\sum_{m\in [M]}\alpha_{nm}= 1,\\
	&&\alpha_{nm},\beta_{nk}\geq 0, \quad \forall m\in [M],k\in[K_n],
\end{eqnarray}  

Accordingly, if only the first case is considered, we can reformulate model PCO as follows.
\begin{subequations} 
	\begin{align}
		\quad & \max_{\boldsymbol{X}\in \mathcal{X}} \frac{1}{N}\sum_{n\in [N]}y_n \\\
		\rm{s.t.} \quad & C(y_n-1)\leq \boldsymbol{\beta}_n^T\boldsymbol{\hat{r}}_n, \quad  \forall n\in [N],\label{W7}\\
		& \boldsymbol{\beta}_n^T\boldsymbol{Q}_n-\boldsymbol{\alpha}^T_n\boldsymbol{X}=0\\
		& \boldsymbol{e}^T\boldsymbol{\alpha}_n= 1,\quad  \forall n\in [N],\\
		& \boldsymbol{\alpha}_n,\boldsymbol{\beta}_n\geq 0,\quad  \forall n\in [N],m\in [M]\\
		& \sum_{l\in\mathcal{L}_a}x_{ml}=1,\forall a\in[A],m\in[M],\\
		&x_{ml},y_n \in \left\{ 0,1 \right\},\quad  \forall n\in [N], l\in[L],m\in[M].
	\end{align}
\end{subequations}

However, there is another case that the minimal $c_n$ is unbounded and therefore, its dual problem W-Utility-D is infeasible. In this case, the infeasibility of model W-Utility-D can only be a consequence of violating constraint \eqref{W4}. Considering the fact that $\inf\limits_{\hat{u}_{nl}\in \mathcal{U}_n(\boldsymbol{Q}_n,\boldsymbol{\hat{r}}_n)}\boldsymbol{1}\{\max_{m\in[M]}\hat{u}_n^T\boldsymbol{x}_m\geq0\}=0$ if the minimal $c_n$ is unbounded and that $\sum_{m\in [M]} x_{ml}\alpha_{nm}\leq 1$,  we can use $y_n$ to relax this infeasibility.
\begin{equation}
	y_n-1\leq \sum_{k\in [K_n]}q_{nlk}\beta_{nk}-\sum_{m\in [M]} x_{ml}\alpha_{nm}\leq 1-y_n, \quad \forall l\in [L].\label{W8}
\end{equation}
When $y_n=0$, $\beta_{nk}=0,\forall k\in[K_n]$ is a feasible solution to the inequality.
Note that we can substitute inequality \eqref{W7} by \eqref{W8} without influence the optimization. This is because the focal situation, $y_n=1$, only indicates the situation where the minimal $c_n$ is bounded.

Through the substitution, we complete the proof. \qed 
\subsection{Proof of Proposition \ref{pro4}}
Because $\bar{u}\in\mathcal{U}_n(\boldsymbol{Q},\boldsymbol{\hat{r}}_n)$, $\boldsymbol{Q}\bar{u}_n=\boldsymbol{\hat{r}}_n$. Therefore, $\boldsymbol{\beta}_n^T\boldsymbol{Q}\bar{u}_n=\boldsymbol{\beta}_n^T\boldsymbol{\hat{r}}_n$. According to equation (\ref{5c}), $\boldsymbol{\beta}_n^T\boldsymbol{Q}=\boldsymbol{x}^T$. As such, $x^T\bar{u}_n=\boldsymbol{\beta}_n^T\boldsymbol{\hat{r}}_n$, which means that when $\boldsymbol{x}$ is feasible to model PCO-RT-S, the ground truth utility is equal to $\boldsymbol{\beta}_n^T\boldsymbol{\hat{r}}_n$. \qed

\begin{lemma}\label{pro:outofsample2}
We have
\begin{itemize}
\item[(i)] With probability higher than $1-\delta$, 
\begin{align*}
 P(\boldsymbol{\tilde{u}}^T\boldsymbol{x^{N,e}})-P(\boldsymbol{\tilde{u}}^T\boldsymbol{x}^*_e) \leq \epsilon 
\end{align*}
when $N\geq \frac{2e^2(L+2)\mbox{log}2}{\epsilon^2}\mbox{log}(\frac{2}{\delta})$
\item [(ii)] \begin{align*}
\mathbb{E}\{ |P(\boldsymbol{\tilde{u}}^T\boldsymbol{x^{N,e}}-P(\boldsymbol{\tilde{u}}^T\boldsymbol{x}^*_e)| \}\leq 2e\mbox{log2}\sqrt{2(L+2)/N}.
\end{align*}
\end{itemize}
\end{lemma}

\subsection{Proof of Lemma \ref{pro:outofsample2}}
Based on Proposition \ref{pro4}, PCO-RT-S model reduces to
\begin{align*}
	\boldsymbol{x}^{N,e}=\mbox{argmax}_{\boldsymbol{x}\in \mathcal{X}^e}  \frac{1}{N}\sum_{n=1}^N \mathbb{I}(\hat{\boldsymbol{u}}_n^T\boldsymbol{x}\geq 0).
\end{align*}
Given a fixed $\boldsymbol{x}\in\mathcal{X}^e$, $\mathbb{I}(\boldsymbol{\hat{u}}_n^Tx\geq 0)$ has only two possible values, i.e., $0$ and $1$.
Then Hoeffding's inequality implies for any $\epsilon\geq 0$
\begin{align*}
	P(  \left|\frac{1}{N}\sum_{n=1}^N \mathbb{I}(\boldsymbol{\hat{u}}_n^T\boldsymbol{x}\geq 0)-P(\boldsymbol{\tilde{u}}^Tx) \right|\geq \epsilon\theta_1\leq 2e^{-2n\epsilon^2}.
\end{align*}
For any random variable $X$, we may consider its \emph{Orlicz norm} $\| X \|_{\psi}$ defined as 
\begin{align*}
	\| X \|_{\psi}=\inf \{C>0: \mathbb{E}\psi( |X|/C \theta_1\leq 1   \}
\end{align*}
where $\psi$ is a nondecreasing, convex function with $\psi(0)=0$.
In this paper, we focus on a specific $\psi(x)\triangleq e^{x^2}-1$.
The Hoeffding's inequality implies that the conditions of Lemma 2.2.1 in \cite{vandervarrt1996} holds with $K=p=2$, $C=2n$,  and $X=\frac{1}{N}\sum_{n=1}^N \mathbb{I}(\boldsymbol{\hat{u}}_n^T\boldsymbol{x}\geq 0)-P(\boldsymbol{\tilde{u}}^Tx)$ with their notations.
And so $\| \frac{1}{N}\sum_{n=1}^N \mathbb{I}(\boldsymbol{\hat{u}}_n^T\boldsymbol{x}\geq 0)-P(\boldsymbol{\tilde{u}}^Tx)\|_{\psi}\leq \sqrt{\frac{3}{2N}} $ for every given $x$ by Lemma 2.2.1 in \cite{vandervarrt1996}.

Next, we verify conditions in Lemma 2.2.2 of  \cite{vandervarrt1996}.
Specifically, for $\psi(x)=e^{x^2}-1$, we claim that
\begin{align*}
	\psi(x)\psi(y)\leq \psi(exy)
\end{align*}
and the proof is omitted for brevity.
Note that with $L$ attributes, the total number of elements in $\boldsymbol{\mathcal{X}^e}$
Then the proofs of Lemma 2.2.2 of  \cite{vandervarrt1996} implies 
\begin{align}
	\| \sup_{\boldsymbol{x}\in \mathcal{X}^e}\left|\frac{1}{N}\sum_{n=1}^N \mathbb{I}(\boldsymbol{\hat{u}}_n^T\boldsymbol{x}\geq 0)-P(\boldsymbol{\tilde{u}}^T\boldsymbol{x})  \right|    \|_{\psi}&\leq e\sqrt{\mbox{log}(1+2^{K+1})}\sup_{x\in\mathcal{X}}\|  \frac{1}{N}\sum_{n=1}^N \mathbb{I}(\boldsymbol{\hat{u}}_n^T\boldsymbol{x}\geq 0)-P(\boldsymbol{\tilde{u}}^Tx)\|_{\psi} \notag\\
	&\leq e\sqrt{3\mbox{log}(1+2^{L+1})/2N}\notag\\
	&\leq e\sqrt{2(L+2)\mbox{log 2}/N}.\label{eq:bounded psi 2 norm}
\end{align}
Denote $x^*=\mbox{argma}_{x\in\boldsymbol{\mathcal{X}^e}} P(\boldsymbol{\tilde{u}}^Tx)$.
Using the inequality that 
\begin{align*}
	\| \max_{x\in\mathcal{X}}g(x)-\max_{x\in\mathcal{X}}\tilde{g}(x) \|\leq \sup_{x\in\mathcal{X}}|g(x)-\tilde{g}(x)  |,
\end{align*}
we conclude for any $\epsilon>0$
\begin{align*}
	P(\left|\frac{1}{N}\sum_{n=1}^N \mathbb{I}(\boldsymbol{\hat{u}}_n^T\boldsymbol{x}^{N,e}\geq 0)-P(\boldsymbol{\tilde{u}}^T\boldsymbol{x}^*)\right| \geq \epsilon\theta_1&=P( \left|\max_{\boldsymbol{x}\in\boldsymbol{\mathcal{X}^e}}\frac{1}{N}\sum_{n=1}^N \mathbb{I}(\boldsymbol{\hat{u}_n}^T\boldsymbol{x}\geq 0)-\max_{\boldsymbol{x}\in\boldsymbol{\mathcal{X}}^e}P(\boldsymbol{\tilde{u}}^Tx^*)    \right|\geq \epsilon \theta_1\\
	&\leq P(  \sup_{\boldsymbol{x}\in \boldsymbol{\mathcal{X}^e}}\left|\frac{1}{N}\sum_{n=1}^N \mathbb{I}(\boldsymbol{\hat{u}}_n^T\boldsymbol{x}\geq 0)-P(\boldsymbol{\tilde{u}}^Tx)  \right| \geq \epsilon \theta_1
\end{align*}
Denote $\Delta_N=\sup_{\boldsymbol{x}\in \boldsymbol{\mathcal{X}^e}}\left|\frac{1}{N}\sum_{n=1}^N \mathbb{I}(\boldsymbol{\hat{u}}_n^T\boldsymbol{x}\geq 0)-P(\boldsymbol{\tilde{u}}^Tx)  \right| $.
Since $\psi$ is nondecreasing, the RHS of the above display can be bounded by
\begin{align*}
	P(  \Delta_N \geq \epsilon \theta_1\leq P(\psi(\Delta_N/\|\Delta_N\|_{\psi})\geq \psi(\epsilon/\|\Delta_N\|_{\psi}))\leq \frac{\mathbb{E}\psi(\Delta_N/\|\Delta_N\|_{\psi}) }{\psi(\epsilon/\|\Delta_N\|_{\psi})}\leq \frac{1}{\psi(\epsilon/\|\Delta_N\|_{\psi})}.
\end{align*}
The last inequality is from the definition of Orlicz norm.
Now that $\|  \Delta_N\|_\psi \leq e\sqrt{2(L+2)\mbox{log 2}/N}$, the last term in the above display can be bounded by 
\begin{align}\label{eq:maximal inequalty}
	P(  \Delta_N \geq \epsilon \theta_1\leq \frac{1}{\psi(\epsilon/\|\Delta_N\|_{\psi})}\leq \frac{1}{e^{\frac{N\epsilon^2}{2e^2(L+2)\mbox{log}2 }}-1}
\end{align}
For $\delta\in (0,1)$, we select $n$ large such that $e^{\frac{N\epsilon^2}{2e^2(K+2)\mbox{log}2 }}\geq 2/\delta $.
Then the above display will be bounded by $\delta$.
This requires an $N\geq \frac{2e^2(L+2)\mbox{log}2}{\epsilon^2}\mbox{log}(\frac{2}{\delta})$.
Now use the telescoping and triangle inequality
\begin{align*}
	|P(\boldsymbol{\tilde{u}}^T\boldsymbol{x^{N,e}}-P(\boldsymbol{\tilde{u}}^T\boldsymbol{x}^*_e)|\leq \left|\frac{1}{N}\sum_{n=1}^N \mathbb{I}(\boldsymbol{\hat{u}}_n^T\boldsymbol{x}^{N,e}\geq 0)-P(\boldsymbol{\tilde{u}}^T\boldsymbol{x}^*_e)  \right|+\left|\frac{1}{N}\sum_{n=1}^N \mathbb{I}(\boldsymbol{\hat{u}}_n^T\boldsymbol{x^{N,e}}\geq 0)-P(\boldsymbol{\tilde{u}}^T\boldsymbol{x^{N,e}})  \right|
\end{align*}
and use \eqref{eq:maximal inequalty} again to conclude.

Now, by the above display and Section 2.2 of \cite{vandervarrt1996}, we know 
\begin{align*}
\mathbb{E}\{ |P(\boldsymbol{\tilde{u}}^T\boldsymbol{x^{N,e}}-P(\boldsymbol{\tilde{u}}^T\boldsymbol{x}^*_e)| \}\leq 2\sqrt{\mbox{log}2}\mathbb{E}\| \sup_{\boldsymbol{x}\in \mathcal{X}^e}\big|\frac{1}{N}\sum_{n=1}^N \mathbb{I}(\boldsymbol{\hat{u}}_n^T\boldsymbol{x}\geq 0)-P(\boldsymbol{\tilde{u}}^T\boldsymbol{x})  \big|    \|_{\psi}.
\end{align*}
Here, the expectation is taken with respect to the randomness in $\boldsymbol{x^{N,e}}$.
Now, use \eqref{eq:bounded psi 2 norm} to see
\begin{align*}
\mathbb{E}\{ |P(\boldsymbol{\tilde{u}}^T\boldsymbol{x^{N,e}}-P(\boldsymbol{\tilde{u}}^T\boldsymbol{x}^*_e)| \}\leq 2e\mbox{log 2}\sqrt{2(L+2)/N}.
\end{align*}
\subsection{Proof Proposition \ref{pro:outofsample1}}
When the distance between any two points in $\mathcal{U}_n(\boldsymbol{Q}_n,\boldsymbol{\hat{r}}_n)$ is no longer than $d, \forall n\in[N]$ and that for every $\boldsymbol{X}$, each point in the uncertainty set is contained in a ball centered at $\boldsymbol{\hat{u}_n} $ with radius $d$, which we denoted as $B(\boldsymbol{\hat{u}_n},d)$.
Then, we can derive
\begin{align*}
&\left|P(\min_{u\in\mathcal{U}_n(\boldsymbol{Q}_n,\boldsymbol{\hat{r}}_n)} \max_{m\in[M]} {{u}}^T\boldsymbol{x}_m\geq 0)-P( \max_{m\in[M]} \boldsymbol{\tilde{u}}^T\boldsymbol{x}_m\geq 0 )\right| \\
&=P( \max_{m\in[M]} \boldsymbol{\tilde{u}}^T\boldsymbol{x}_m\geq 0 )-P(\min_{u\in\mathcal{U}_n(\boldsymbol{Q}_n,\boldsymbol{\hat{r}}_n)} \max_{m\in[M]} {{u}}^T\boldsymbol{x}_m\geq 0) \\
&\leq P(\max_{u\in B(\boldsymbol{\tilde{u}},d)}  \max_{m\in[M]}{{u}}^T\boldsymbol{x}_m\geq 0)-P(\min_{u\in{B}(\boldsymbol{\tilde{u}},d)} \max_{m\in[M]} {{u}}^T\boldsymbol{x}_m\geq 0)\\
&= P(\min_{u\in{B}(\boldsymbol{\tilde{u}},d)} \max_{m\in[M]} {{u}}^T\boldsymbol{x}_m\leq 0,\max_{u\in B(\boldsymbol{\tilde{u}},d)}  \max_{m\in[M]}{{u}}^T\boldsymbol{x}_m\geq 0 )\\
&= 1-P(\min_{u\in{B}(\boldsymbol{\tilde{u}},d)} \max_{m\in[M]} {{u}}^T\boldsymbol{x}_m\geq 0)-P(\max_{u\in B(\boldsymbol{\tilde{u}},d)}  \max_{m\in[M]}{{u}}^T\boldsymbol{x}_m\leq  0)
\end{align*}
The first equality holds due to the fact $P(\min_{u\in\mathcal{U}_n(\boldsymbol{Q}_n,\boldsymbol{\hat{r}}_n)} \max_{m\in[M]} {{u}}^T\boldsymbol{x}_m\geq 0)-P( \max_{m\in[M]} \boldsymbol{\tilde{u}}^T\boldsymbol{x}_m\geq 0 )\leq 0 $; the second inequality holds because $P( \max_{m\in[M]} \boldsymbol{u}^T\boldsymbol{x}_m\geq 0 )\leq P(\max_{u\in \mathcal{U}_n(\boldsymbol{Q}_n,\boldsymbol{\hat{r}}_n)} \max_{m\in[M]} \boldsymbol{u}^T\boldsymbol{x}_m\geq 0) $; the third equality holds because $\{\tilde{u}|\min_{u\in \mathcal{U}_n(\boldsymbol{Q}_n,\boldsymbol{\hat{r}}_n)} \max_{m\in[M]} \boldsymbol{u}^T\boldsymbol{x}_m\geq 0\}\subseteq \{\tilde{u}|\max_{u\in B(\boldsymbol{\tilde{u}},d)} \max_{m\in[M]} \boldsymbol{u}^T\boldsymbol{x}_m\geq 0\}$; and the last equality holds because the three events are disjoint and compose the sample space of $\tilde{u}$.
Then, we focus on the problem
\begin{eqnarray}
	&\max_{u}& \max_{m\in[M]} \boldsymbol{u}^T\boldsymbol{x}_m\\
	&s.t.& \sum_{l\in[L]}(u_l-\tilde{u}_l)^2\leq d^2,
\end{eqnarray}
which is equivalent to
\begin{eqnarray}
	&\max_{\delta}& \max_{m\in[M]} (\boldsymbol{\tilde{u}}^T+\boldsymbol{\delta}^T)\boldsymbol{x}_m\\
	&s.t.& \sum_{l\in[L]}\delta_l^2\leq d^2
\end{eqnarray}
Because $\max_{m\in[M]} (\boldsymbol{\tilde{u}}^T+\boldsymbol{\delta}^T)\boldsymbol{x}_m\leq  \max_{m\in[M]} \boldsymbol{\tilde{u}}^T\boldsymbol{x}_m+\max_{m\in[M]}\boldsymbol{\delta}^T\boldsymbol{x}_m$, the above problem is upper bounded by 
\begin{eqnarray}
	&\max_{\delta}& \max_{m\in[M]} \boldsymbol{\tilde{u}}^T\boldsymbol{x}_m+\max_{m\in[M]}\boldsymbol{\delta}^T\boldsymbol{x}_m\\
	&s.t.& \sum_{l\in[L]}\delta_l^2\leq d^2
\end{eqnarray}
The first term is independent of $\delta$. As such, it is equivalent to the following problem.
\begin{eqnarray}
	&\max_{\delta}& \max_{m\in[M]}\boldsymbol{\delta}^T\boldsymbol{x}_m\\
	&s.t.& \sum_{l\in[L]}\delta_l^2\leq d^2
\end{eqnarray}
We can solve the problem individually for all $m\in[M]$ and then find the maximal objective over all $m$.
\begin{eqnarray}
	&\max_{\delta_m}&\boldsymbol{\delta_m}^T\boldsymbol{x}_m\\
	&s.t.& \sum_{l\in[L]}\delta_{lm}^2\leq d^2,\forall m\in[M].
\end{eqnarray}
It is easy to prove the optimal objective value for the above problem is $\sqrt{A}d$ when $\sum_{l\in \mathcal{L}_a}x_{ml}=1,\forall m\in[M],a\in[A]$. As such, $\max_{u} \max_{m\in[M]} \boldsymbol{\tilde{u}}^T\boldsymbol{x}_m\leq \max_{m\in[M]} \boldsymbol{\tilde{u}}^T\boldsymbol{x}_m+\sqrt{A}d$, and accordingly,
$$
P( \max_{m\in[M]} \boldsymbol{\tilde{u}}^T\boldsymbol{x}_m \leq -\sqrt{A}d )\leq P(\max_{u\in B(\boldsymbol{\tilde{u}},d)}  \max_{m\in[M]}{{u}}^T\boldsymbol{x}_m\leq 0).
$$

Similarly, it is easy to prove that
$$\min_{u\in{B}(\boldsymbol{\tilde{u}},d)} \max_{m\in[M]} {{u}}^T\boldsymbol{x}_m\geq \max_{m\in[M]} \min_{u_m\in{B}(\boldsymbol{\tilde{u}},d)}  {{u}}_m^T\boldsymbol{x}_m=\max_{m\in[M]} \boldsymbol{\tilde{u}}^T\boldsymbol{x}_m-\sqrt{A}d.$$
and 
$$
P( \max_{m\in[M]} \boldsymbol{\tilde{u}}^T\boldsymbol{x}_m\geq \sqrt{A}d )\leq P(\min_{u\in{B}(\boldsymbol{\tilde{u}},d)} \max_{m\in[M]} {{u}}^T\boldsymbol{x}_m\geq 0)
$$
Therefore, 
\begin{align*}
&1-P(\min_{u\in{B}(\boldsymbol{\tilde{u}},d)} \max_{m\in[M]} {{u}}^T\boldsymbol{x}_m\geq 0)-P(\max_{u\in B(\boldsymbol{\tilde{u}},d)}  \max_{m\in[M]}{{u}}^T\boldsymbol{x}_m\leq 0)\\
&\leq 1-P( \max_{m\in[M]} \boldsymbol{\tilde{u}}^T\boldsymbol{x}_m\geq \sqrt{A}d )-P( \max_{m\in[M]} \boldsymbol{\tilde{u}}^T\boldsymbol{x}_m\leq -\sqrt{A}d )\\
&=P(-\sqrt{A}d \leq   \max_{m\in[M]} \boldsymbol{\tilde{u}}^T\boldsymbol{x}_m\leq \sqrt{A}d ).
\end{align*}
Suppose the CDF of $\max_{m\in[M]} \boldsymbol{\tilde{u}}^T\boldsymbol{x}_m$ is Lipschitz with a common Lipschitz constant $\theta$ for every $\boldsymbol{x}\in\boldsymbol{X}$.
We have
\begin{align}\label{eq:bounded diff poly}
\left|P(\min_{u\in \mathcal{U}_n(\boldsymbol{Q}_n,\boldsymbol{\hat{r}}_n)} \max_{m\in[M]} {{u}}^T\boldsymbol{x}_m\geq 0)-P( \max_{m\in[M]} \boldsymbol{\tilde{u}}^T\boldsymbol{x}_m\geq 0 )\right|\leq \theta\sqrt{A}d,\forall \boldsymbol{X}\in\mathcal{X}.
\end{align}
Now, let ${\boldsymbol{x}}^N_m\triangleq \mbox{argmax}_{\boldsymbol{x}}\frac{1}{N}\sum_{n=1}^N \mathbb{I}(\min_{u\in \mathcal{U}_n(\boldsymbol{Q}_n,\boldsymbol{\hat{r}}_n)} \max_{m\in[M]} {{u}}^T\boldsymbol{x}_m\geq 0) $ and $$\underline{\boldsymbol{x}}_m^*\triangleq \mbox{argmax}_{\boldsymbol{x}} P(\min_{u\in \mathcal{U}_n(\boldsymbol{Q}_n,\boldsymbol{\hat{r}}_n)} \max_{m\in[M]} {{u}}^T\boldsymbol{x}_m\geq 0).$$
Similar to the proof of Lemma \ref{pro:outofsample2} before, we have
\begin{align*}
P(\left|\frac{1}{N}\sum_{n=1}^N \mathbb{I}(\min_{u\in \mathcal{U}_n(\boldsymbol{Q}_n,\boldsymbol{\hat{r}}_n)} \max_{m\in[M]} {{u}}^T{{\boldsymbol{x}}^N_m}\geq 0)  -P(\min_{u\in \mathcal{U}_n(\boldsymbol{Q}_n,\boldsymbol{\hat{r}}_n)} \max_{m\in[M]} {{u}}^T\underline{\boldsymbol{x}}_m^*\geq 0)\right|\geq \epsilon\theta_1\leq 1-\delta
\end{align*}
for $N\geq \frac{2e^2(LM+2)\mbox{log}(2)}{\epsilon^2}\mbox{log}(\frac{2}{\delta})$.
On the other hand, 
let $\overline{\boldsymbol{x}}_m^N\triangleq \mbox{argmax}_{\boldsymbol{x}}\frac{1}{N}\mathbb{I}(\sum_{n=1}^N  \max_{m\in[M]} {{u}}^T\boldsymbol{x}_m\geq 0) .$
We have
\begin{align}\label{eq:ref 1}
P(\left|\frac{1}{N}\sum_{n=1}^N \mathbb{I}( \max_{m\in[M]} \boldsymbol{\hat{u}}_n^T\overline{\boldsymbol{x}}^N_m\geq 0)  -P(\max_{i=1,\cdots,M} \boldsymbol{\tilde{u}}^T{{\boldsymbol{x}}_m^*}\geq 0)\right|\geq \epsilon\theta_1\leq 1-\delta
\end{align}
for $N\geq \frac{2e^2(LM+2)\mbox{log}(2)}{\epsilon^2}\mbox{log}(\frac{2}{\delta})$.
Write $P( \max_{m\in[M]} {\boldsymbol{\tilde{u}}}^T{{\boldsymbol{x}}_m^N}\geq 0)-P( \max_{m\in[M]} {\boldsymbol{\tilde{u}}}^T{{\boldsymbol{x}}_m^*}\geq 0)$ as
\begin{align}
&P( \max_{m\in[M]} {\boldsymbol{\tilde{u}}}^T{{\boldsymbol{x}}_m^N}\geq 0)-\frac{1}{N}\sum_{n=1}^N \mathbb{I}( \max_{m\in[M]} \boldsymbol{\hat{u}}_n^T{\boldsymbol{x}}^N_m\geq 0)+\notag\\
&\frac{1}{N}\sum_{n=1}^N \mathbb{I}( \max_{m\in[M]} \boldsymbol{\hat{u}}_n^T{\boldsymbol{x}}^N_m\geq 0)-\frac{1}{N}\sum_{n=1}^N \mathbb{I}( \max_{m\in[M]} \boldsymbol{\hat{u}}_n^T\bar{\boldsymbol{x}}^N_m\geq 0)+\notag\\
&\frac{1}{N}\sum_{n=1}^N \mathbb{I}( \max_{m\in[M]} \boldsymbol{\hat{u}}_n^T\bar{\boldsymbol{x}}^N_m\geq 0)-P(\max_{i=1,\cdots,M} \boldsymbol{\tilde{u}}^T{{\boldsymbol{x}}_m^*}\geq 0).\label{eq:decompose_main}
\end{align} 
The first term can be similarly bounded as Lemma \ref{pro:outofsample2} before and the last term has been investigated above in \eqref{eq:ref 1}. 
Thus, with high probability the above display satisfy
\begin{align}\label{eq:decompose_first}
\frac{1}{N}\sum_{n=1}^N \mathbb{I}( \max_{m\in[M]} \boldsymbol{\hat{u}}_n^T{\boldsymbol{x}}^N_m\geq 0)-\frac{1}{N}\sum_{n=1}^N \mathbb{I}( \max_{m\in[M]} \boldsymbol{\hat{u}}_n^T\bar{\boldsymbol{x}}^N_m\geq 0)-2\epsilon\leq P( \max_{m\in[M]} {\boldsymbol{\tilde{u}}}^T{{\boldsymbol{x}}_m^N}\geq 0)-P( \max_{m\in[M]} {\boldsymbol{\tilde{u}}}^T{{\boldsymbol{x}}_m^*}\geq 0)\leq 0
\end{align}
For the difference on the LHS, it is non-positive by definition of  $\bar{x}_m^N$, we further write it as
\begin{align}\label{eq:decompose_second}
&\frac{1}{N}\sum_{n=1}^N \mathbb{I}( \max_{m\in[M]} \boldsymbol{\hat{u}}_n^T{\boldsymbol{x}}^N_m\geq 0)-\frac{1}{N}\sum_{n=1}^N \mathbb{I}(\min_{u\in \mathcal{U}_n(\boldsymbol{Q}_n,\boldsymbol{\hat{r}}_n)} \max_{m\in[M]} {{u}}^T{{\boldsymbol{x}}^N_m}\geq 0)\notag\\
&+\frac{1}{N}\sum_{n=1}^N \mathbb{I}(\min_{u\in \mathcal{U}_n(\boldsymbol{Q}_n,\boldsymbol{\hat{r}}_n)} \max_{m\in[M]} {{u}}^T{{\boldsymbol{x}}^N_m}\geq 0)-P(\min_{u\in \mathcal{U}_n(\boldsymbol{Q}_n,\boldsymbol{\hat{r}}_n)} \max_{m\in[M]} {{u}}^T{\underline{\boldsymbol{x}}^*_m}\geq 0)\notag\\
&+P(\min_{u\in \mathcal{U}_n(\boldsymbol{Q}_n,\boldsymbol{\hat{r}}_n)} \max_{m\in[M]} {{u}}^T{\underline{\boldsymbol{x}}^*_m}\geq 0)-P( \max_{m\in[M]} \boldsymbol{\tilde{u}}^T{\boldsymbol{x}}^*_i\geq 0)\notag\\
&+P( \max_{m\in[M]} \boldsymbol{\tilde{u}}^T{\boldsymbol{x}}^*_i\geq 0)-\frac{1}{N}\sum_{n=1}^N \mathbb{I}( \max_{m\in[M]} \boldsymbol{\hat{u}}_n^T\bar{\boldsymbol{x}}^N_m\geq 0).
\end{align} 
The first line is nonnegative by definition of $\mathcal{U}_n(\boldsymbol{Q}_n,\boldsymbol{\hat{r}}_n) $.
The second and fourth lines are both differences between SAA and true mean.
When $N\geq \frac{2e^2(LM+2)\mbox{log}(2)}{\epsilon^2}\mbox{log}(\frac{2}{\delta})$, the probability that these two lines are above $-2\epsilon$ is higher than $1-2\delta$. 
The third line is deterministic and can be controlled by \eqref{eq:bounded diff poly}.
So with high probability the whole display is bounded below by $-2\epsilon-\theta\sqrt{A}d$.
Now from \eqref{eq:decompose_first}, we know with high probability
\begin{align*}
\theta\sqrt{A}d-4\epsilon\leq P( \max_{m\in[M]} {\boldsymbol{\tilde{u}}}^T{{\boldsymbol{x}}_m^N}\geq 0)-P( \max_{m\in[M]} {\boldsymbol{\tilde{u}}}^T{{\boldsymbol{x}}_m^*}\geq 0)\leq 0
\end{align*}

Combine all and use union bound, we have for $N\geq \frac{32e^2(LM+2)\mbox{log}(2)}{\epsilon^2}\mbox{log}(\frac{4}{\delta})$,
\begin{align*}
P(\left| P( \max_{m\in[M]} {\boldsymbol{\tilde{u}}}^T{{\boldsymbol{x}}_m^N}\geq 0)-P( \max_{m\in[M]} {\boldsymbol{\tilde{u}}}^T{{\boldsymbol{x}}_m^*}\geq 0) \right|>\epsilon+\theta\sqrt{A}d\theta_1<1-\delta.
\end{align*}

Now, from the fact $|\sup f-\sup g|\leq \sup|f-g|$ and \eqref{eq:decompose_main}, we have
\begin{align}
\mathbb{E}\left| P( \max_{m\in[M]} {\boldsymbol{\tilde{u}}}^T{{\boldsymbol{x}}_m^N}\geq 0)-P( \max_{m\in[M]} {\boldsymbol{\tilde{u}}}^T{{\boldsymbol{x}}_m^*}\geq 0) \right|\leq & 2\mathbb{E} \sup_x\left |\frac{1}{N}\sum_{n=1}^N \mathbb{I}( \max_{m\in[M]} \boldsymbol{\hat{u}}_n^T\bar{\boldsymbol{x}}_m\geq 0)-P(\max_{i=1,\cdots,M} \boldsymbol{\tilde{u}}^T{{\boldsymbol{x}}_m}\geq 0)  \right|\notag\\
&+\mathbb{E}\left|\frac{1}{N}\sum_{n=1}^N \mathbb{I}( \max_{m\in[M]} \boldsymbol{\hat{u}}_n^T{\boldsymbol{x}}^N_m\geq 0)-\frac{1}{N}\sum_{n=1}^N \mathbb{I}( \max_{m\in[M]} \boldsymbol{\hat{u}}_n^T\bar{\boldsymbol{x}}^N_m\geq 0)  \right|\label{eq:decompose expectation err}
\end{align}
The first term on the RHS can be analysed similar to (ii) of Lemma \ref{pro:outofsample2}. 
Next, the second term on the RHS can be analysed by \eqref{eq:decompose_second}.
By \eqref{eq:decompose_second} and Section 2.2 of \cite{vandervarrt1996}, we would have
\begin{align*}
&\mathbb{E}\left|\frac{1}{N}\sum_{n=1}^N \mathbb{I}( \max_{m\in[M]} \boldsymbol{\hat{u}}_n^T{\boldsymbol{x}}^N_m\geq 0)-\frac{1}{N}\sum_{n=1}^N \mathbb{I}( \max_{m\in[M]} \boldsymbol{\hat{u}}_n^T\bar{\boldsymbol{x}}^N_m\geq 0)  \right|\leq \theta \sqrt{A} d\\
&+\sqrt{\mbox{log}2}\mathbb{E}\| \frac{1}{N}\sum_{n=1}^N \mathbb{I}(\min_{u\in \mathcal{U}_n(\boldsymbol{Q}_n,\boldsymbol{\hat{r}}_n)} \max_{m\in[M]} {{u}}^T{{\boldsymbol{x}}^N_m}\geq 0)-P(\min_{u\in \mathcal{U}_n(\boldsymbol{Q}_n,\boldsymbol{\hat{r}}_n)} \max_{m\in[M]} {{u}}^T{\underline{\boldsymbol{x}}^*_m}\geq 0) \|_{\psi}\\
&+\sqrt{\mbox{log}2}\mathbb{E}\| P( \max_{m\in[M]} \boldsymbol{\tilde{u}}^T{\boldsymbol{x}}^*_i\geq 0)-\frac{1}{N}\sum_{n=1}^N \mathbb{I}( \max_{m\in[M]} \boldsymbol{\hat{u}}_n^T\bar{\boldsymbol{x}}^N_m\geq 0 \|_{\psi}
\end{align*}
The terms $$\mathbb{E} \sup_x\left |\frac{1}{N}\sum_{n=1}^N \mathbb{I}( \max_{m\in[M]} \boldsymbol{\hat{u}}_n^T\bar{\boldsymbol{x}}_m\geq 0)-P(\max_{i=1,\cdots,M} \boldsymbol{\tilde{u}}^T{{\boldsymbol{x}}_m}\geq 0)  \right|$$ and  $$\mathbb{E}\| \frac{1}{N}\sum_{n=1}^N \mathbb{I}(\min_{u\in \mathcal{U}_n(\boldsymbol{Q}_n,\boldsymbol{\hat{r}}_n)} \max_{m\in[M]} {{u}}^T{{\boldsymbol{x}}^N_m}\geq 0)-P(\min_{u\in \mathcal{U}_n(\boldsymbol{Q}_n,\boldsymbol{\hat{r}}_n)} \max_{m\in[M]} {{u}}^T{\underline{\boldsymbol{x}}^*_m}\geq 0) \|_{\psi} $$ both characterize the difference between SAA and true mean, and can both be bounded following  same logic in the proof Lemma \ref{pro:outofsample2}.
Then in view of \eqref{eq:decompose expectation err}, we conclude
\begin{align*}
\mathbb{E}\left| P( \max_{m\in[M]} {\boldsymbol{\tilde{u}}}^T{{\boldsymbol{x}}_m^N}\geq 0)-P( \max_{m\in[M]} {\boldsymbol{\tilde{u}}}^T{{\boldsymbol{x}}_m^*}\geq 0) \right|\leq \theta \sqrt{A} d+4{\mbox{log}2} \epsilon\sqrt{2(LM+2)/N}.
\end{align*}

.

\section{Mean Vector and Covariance Matrix of the Distribution of Part-worth Utilities}\label{appendix_data}
The mean vector and covariance matrix of part-worth utilities' the distribution presented in \cite{Gilbride2004} are shown in Tables 
\ref{tab Mean of partworth utilities} and \ref{tab Covariance}.
\begin{table}[htbp]
\caption{Mean of part-worth utilities}
\centering
\scriptsize
\label{tab Mean of partworth utilities}
\begin{tabular}{llr}
	\hline
	Attributes (a)& Levels (l) & Mean($\bar{u}_l$) \\
	\hline
	Body Sytle (1)& Low body (1)& 0.325\\
	& Medium body (2) & 2.843\\
	& High body (3) & 2.224\\
	Midroll change (2) & Manual change (4)& 0.059\\
	& Automatic change (5) & 0.178\\
	Annotation (3)& Preset list (6) & 0.362\\
	& Customerized list (7)& 0.678\\
	& Input method 1 (8) & -1.036\\
	& Input method 2 (9) & 0.767\\
	& Input method 3 (10)& -0.985\\
	Camera operations feedback (4) & Operation feedback (11) & 0.375\\
	Zoom (5)& 2$\times$ zoom (12)& 1.122\\
	& 4$\times$ zoom (13)& 1.477\\
	Viewfinder (6) & Large viewfinder (14)& -0.158\\
	Setting feedback (7)& LCD (15) & 0.759\\
	& Viewfinder (16) & 0.757\\
	& LCD \& viewfinder (17)& 0.902\\
	\hline 
\end{tabular}
\end{table}

\begin{table}[htbp]
\caption{Covariance matrix of part-worth utility}
\label{tab Covariance}
\tiny
\centering
\begin{tabular}{lrrrrrrrrrrrrrrrrrr}
	\hline
	& 1 & 2 & 3 & 4 & 5& 6 & 7 & 8& 9& 10& 11& 12& 13 & 14 & 15& 16& 17 & \\
	\hline
	1& 19.605 &&&&&&&&&&&&&&&&&\\
	2& 17.952 & 20.258 &&&&&&&&&&&&&&&&\\
	3& 15.303 & 17.674 & 17.845 &&&&&&&&&&&&&&&\\
	4& -0.548 & -0.906 & -1.218 & 1.803&&&&&&&&&&&&&&\\
	5& -2.077 & -2.232 & -2.153 & 0.507& 1.4&&&&&&&&&&&&&\\
	6& -0.575 & -0.309 & -0.333 & 0.097& 0.171 & 0.629&&&&&&&&&&&&\\
	7& -0.584 & -0.366 & -0.356 & 0.147& 0.231 & 0.255& 0.685&&&&&&&&&&&\\
	8& -8.252 & -7.975 & -6.735 & 0.371& 1.134 & 0.3 & 0.387& 5.392 &&&&&&&&&&\\
	9& -1.921 & -1.52& -1.285 & 0.139& 0.395 & 0.237& 0.26& 0.974 & 1.096 &&&&&&&&&\\
	10 & -5.553 & -5.628 & -4.753 & 0.432& 0.854 & 0.238& 0.224& 2.734 & 0.726 & 3.315&&&&&&&&\\
	11 & -0.808 & -0.67& -0.665 & 0.061& 0.173 & 0.019& 0.033& 0.494 & 0.225 & 0.347& 1.055&&&&&&&\\
	12 & -0.207 & 0.651& 0.895& -0.009 & 0.022 & 0.106& 0.121& 0.134 & 0.167 & -0.056 & 0.115& 1.42&&&&&&\\
	13 & -0.811 & 0.54& 1.08& -0.225 & -0.08 & 0.225& 0.173& 0.418 & 0.251 & -0.021 & -0.008 & 1.07& 2.21&&&&&\\
	14 & -1.394 & -1.054 & -0.859 & -0.011 & 0.149 & 0.051& -0.021 & 0.592 & 0.231 & 0.487& 0.064& 0.122& 0.13& 0.947 &&&&\\
	15 & -0.434 & 0.06& 0.324& -0.075 & 0.04& 0.106& 0.12& 0.167 & 0.181 & 0.04& 0.092& 0.343& 0.454 & 0.141 & 0.966&&&\\
	16 & -0.289 & 0.035& 0.161& 0.011& 0.044 & 0.083& 0.067& 0.129 & 0.179 & 0.083& 0.123& 0.298& 0.336 & 0.119 & 0.525& 0.976&&\\
	17 & -0.472 & -0.042 & 0.12& 0.019& 0.057 & 0.109& 0.118& 0.255 & 0.173 & 0.084& 0.106& 0.321& 0.427 & 0.074 & 0.526& 0.538& 0.937&\\
	\hline
\end{tabular}
\end{table}

\section{Orthogonal Design applied in Section \ref{EPCO}}\label{OrthCA}
\begin{table}[hptb]
\caption{The fixed orthogonal designed CA survey}

\centering
\scriptsize
\begin{tabular}{llllllllllllllllllllllllllllllll}
	\hline
	K  & 1  & 2  & 3  & 4  & 5  & 6  & 7  & 8  & 9  & 10 & 11 & 12 & 13 & 14 & 15 & 16 & 17 & 18 & 19 & 20 & 21 & 22 & 23 & 24 & 25 & 26 & 27 & 28 & 29 & 30 & 31 \\ \hline
	1  & 1  & 0  & 0  & -1 & 0  & 0  & 0  & 0  & 0  & 0  & 0  & 0  & 0  & 1  & 0  & 0  & -1 & 0  & 0  & 0  & 0  & 0  & 1  & 0  & 0  & -1 & 0  & 0  & 0  & 0  & 0  \\
	2  & -1 & 0  & 0  & 1  & 0  & 1  & -1 & 0  & 0  & 0  & 0  & 0  & 0  & 0  & 0  & 0  & 1  & -1 & 0  & 0  & 0  & 0  & 0  & 0  & 0  & 1  & 0  & 0  & 0  & -1 & 0  \\
	3  & 1  & 0  & 0  & 0  & -1 & 0  & 1  & 0  & 0  & 0  & 0  & -1 & 0  & 0  & 0  & 0  & 0  & 1  & 0  & 0  & -1 & 0  & 0  & 0  & 0  & 0  & 0  & 0  & 0  & 1  & -1 \\
	4  & 0  & 0  & 0  & -1 & 1  & -1 & 0  & 0  & 0  & 0  & 0  & 1  & 0  & 0  & 0  & 0  & 0  & 0  & 0  & 0  & 0  & 0  & 0  & 0  & 0  & 0  & 0  & -1 & 0  & 0  & 1  \\
	5  & 0  & 0  & 0  & 0  & 0  & 1  & 0  & 0  & 0  & 0  & -1 & 0  & 0  & 0  & -1 & 0  & 0  & 0  & 0  & 0  & 1  & 0  & 0  & 0  & 0  & 0  & 0  & 1  & -1 & 0  & 0  \\
	6 & 0  & 0  & -1 & 1  & 0  & 0  & 0  & -1 & 0  & 0  & 1  & 0  & 0  & 0  & 1  & 0  & 0  & 0  & 0  & 0  & 0  & -1 & 0  & 0  & 0  & 0  & 0  & -1 & 1  & 0  & 0  \\
	7 & 0  & -1 & 1  & 0  & 0  & 0  & 0  & 1  & 0  & 0  & 0  & -1 & 0  & -1 & 0  & 0  & 0  & 0  & 0  & 0  & 0  & 1  & 0  & -1 & 0  & 0  & 0  & 1  & 0  & 0  & 0  \\
	8 & 0  & 0  & 0  & 0  & 0  & 0  & -1 & 0  & 0  & 0  & 0  & 1  & 0  & 0  & 0  & 0  & 0  & 0  & 0  & 0  & 0  & 0  & 0  & 1  & 0  & -1 & 0  & 0  & 0  & 0  & 0  \\
	9 & 0  & 1  & -1 & 0  & 0  & 0  & 0  & 0  & 0  & 0  & 0  & 0  & 0  & 1  & 0  & 0  & 0  & 0  & 0  & -1 & 0  & 0  & -1 & 0  & 0  & 1  & 0  & 0  & 0  & 0  & 0  \\
	10 & 0  & 0  & 1  & -1 & 0  & 0  & 1  & 0  & 0  & 0  & 0  & -1 & 0  & 0  & 0  & -1 & 0  & 0  & 0  & 1  & 0  & 0  & 0  & 0  & 0  & 0  & 0  & 0  & 0  & 0  & 0  \\
	11 & 0  & 0  & 0  & 0  & 0  & 0  & 0  & -1 & 0  & 0  & 0  & 1  & 0  & 0  & 0  & 1  & 0  & 0  & 0  & -1 & 0  & 0  & 1  & 0  & 0  & -1 & 0  & 0  & 0  & 0  & 0  \\
	12 & -1 & 0  & 0  & 1  & 0  & 0  & 0  & 1  & -1 & 0  & 0  & 0  & 0  & 0  & 0  & 0  & 0  & 0  & 0  & 0  & 0  & 0  & 0  & 0  & 0  & 1  & 0  & 0  & 0  & 0  & -1 \\
	13 & 0  & 0  & 0  & 0  & 0  & 0  & 0  & 0  & 1  & 0  & 0  & -1 & 0  & 0  & 0  & 0  & -1 & 0  & 0  & 1  & 0  & 0  & 0  & 0  & 0  & 0  & -1 & 0  & 0  & 0  & 1  \\
	14 & 0  & 0  & 0  & 0  & 0  & -1 & 0  & 0  & 0  & 0  & 0  & 1  & 0  & 0  & 0  & 0  & 1  & 0  & 0  & -1 & 0  & 0  & 0  & 0  & 0  & 0  & 0  & 0  & 0  & 0  & 0  \\
	15 & 1  & 0  & 0  & -1 & 0  & 1  & 0  & 0  & -1 & 0  & 0  & 0  & 0  & -1 & 0  & 0  & 0  & 0  & 0  & 1  & 0  & 0  & 0  & 0  & -1 & 0  & 1  & 0  & 0  & 0  & 0  \\
	16 & 0  & 0  & -1 & 1  & 0  & 0  & 0  & 0  & 1  & -1 & 0  & 0  & 0  & 1  & 0  & 0  & 0  & 0  & 0  & 0  & -1 & 0  & -1 & 0  & 1  & 0  & 0  & 0  & 0  & 0  & 0  \\
	17 & 0  & 0  & 0  & 0  & 0  & 0  & 0  & -1 & 0  & 1  & 0  & 0  & 0  & -1 & 0  & 0  & 0  & 0  & 0  & 0  & 1  & 0  & 1  & 0  & 0  & 0  & 0  & 0  & -1 & 0  & 0  \\
	18 & 0  & -1 & 1  & 0  & 0  & 0  & 0  & 1  & 0  & 0  & 0  & 0  & -1 & 1  & 0  & 0  & -1 & 0  & 0  & 0  & 0  & 0  & 0  & 0  & 0  & 0  & 0  & 0  & 1  & -1 & 0  \\
	19 & 0  & 1  & -1 & 0  & 0  & 0  & -1 & 0  & 0  & 0  & 0  & 0  & 1  & 0  & -1 & 0  & 1  & 0  & 0  & 0  & 0  & 0  & 0  & 0  & 0  & 0  & -1 & 0  & 0  & 1  & 0  \\
	20  & 0  & 0  & 0  & 0  & 0  & 0  & 0  & 1  & -1 & 0  & 0  & 0  & 0  & 0  & 0  & -1 & 1  & 0  & 0  & 0  & 0  & 0  & 1  & -1 & 0  & 0  & 0  & 0  & 0  & 0  & 0  \\
	21  & 0  & 0  & 0  & 0  & 0  & -1 & 1  & 0  & 0  & 0  & 0  & 0  & 0  & 0  & 0  & 0  & 0  & 0  & 1  & 0  & 0  & -1 & 0  & 0  & 0  & 0  & 0  & 0  & 0  & 0  & 0  \\
	22  & 1  & -1 & 0  & 0  & 0  & -1 & 0  & 0  & 0  & 0  & 0  & 0  & 1  & 0  & 0  & 0  & 0  & 0  & 0  & 0  & 0  & 0  & -1 & 0  & 0  & 0  & 0  & 0  & 1  & 0  & 0  \\
	23  & 1  & 0  & 0  & -1 & 0  & 0  & 0  & 0  & 0  & 0  & 0  & 0  & 0  & 0  & 0  & 0  & -1 & 0  & 0  & 0  & 1  & 0  & 0  & 0  & -1 & 1  & 0  & 0  & 0  & 0  & 0  \\
	24  & 0  & 0  & 1  & -1 & 0  & 0  & 0  & 0  & 0  & -1 & 0  & 0  & 1  & 0  & 0  & -1 & 0  & 1  & 0  & 0  & 0  & 0  & 0  & 0  & 0  & 0  & -1 & 0  & 0  & 0  & 1  \\
	25 & 0  & 1  & -1 & 0  & 0  & 0  & 0  & 0  & -1 & 0  & 0  & 0  & 1  & 0  & 0  & 1  & 0  & -1 & 0  & 0  & 0  & 0  & 0  & 0  & 0  & 0  & 0  & 0  & 0  & 0  & 0  \\
	26 & -1 & 0  & 0  & 1  & 0  & 0  & 0  & -1 & 0  & 0  & 0  & 1  & 0  & 0  & -1 & 0  & 0  & 0  & 0  & 1  & 0  & 0  & 0  & 0  & 0  & 0  & 0  & 0  & 0  & 0  & 0  \\
	27 & 0  & 0  & 1  & 0  & -1 & 1  & 0  & 0  & 0  & 0  & -1 & 0  & 0  & 0  & 0  & 0  & 0  & 0  & 0  & 0  & 0  & 0  & 0  & 0  & 1  & 0  & 0  & 0  & 0  & -1 & 0  \\ \hline
\end{tabular}

\end{table}
%%%%%%%%%%%%%%%%%
\end{document}